\documentclass{article}
\usepackage{amssymb}
\usepackage{amsfonts}
\usepackage{amsmath}

\input{tcilatex}
\begin{document}

\title{Metrizability of the Lie algebroid generalized tangent bundle and
(generalized) Lagrange $\left( \rho ,\eta \right) $-spaces}
\author{Constantin M. ARCU\c{S} \\
%EndAName
\ \ \\
[0pt] \ \ \\
[0pt]
\begin{tabular}{c}
SECONDARY SCHOOL \textquotedblleft CORNELIUS RADU\textquotedblright , \\
RADINESTI VILLAGE, 217196, GORJ COUNTY, ROMANIA \\
e-mail: c\_arcus@radinesti.ro, constantin\_arcus@yahoo.ro,%
\end{tabular}%
}
\maketitle

\begin{abstract}
A class of metrizable vector bundles in the general framework of
generalized Lie algebroids have been presented in the eight
reference. Using a generalized Lie algebroid we obtain the Lie
algebroid generalized tangent bundle of a vector bundle. This Lie
algebroid is a new example of metrizable vector bundle. A new class
of Lagrange spaces, called by use, generalized Lagrange $\left( \rho
,\eta \right) $-space, Lagrange $\left( \rho ,\eta \right) $-space
and Finsler $\left( \rho ,\eta \right) $-space are presented. In the
particular case of Lie algebroids, new and important results are
presented. In particular, if all morphisms are identities morphisms,
then the classical results are obtained. \ \ \bigskip\newline
\end{abstract}

\tableofcontents

\section{Introduction}

\ \

In general, if $\mathcal{C}$ is a category, then we denote $\left\vert
\mathcal{C}\right\vert $ the class of objects and for any $A,B{\in }%
\left\vert \mathcal{C}\right\vert $, we denote $\mathcal{C}\left( A,B\right)
$ the set of morphisms of $A$ source and $B$ target and $Iso_{\mathcal{C}%
}\left( A,B\right) $ the set of $\mathcal{C}$-isomorphisms of $A$ source and
$B$ target. Let $\mathbf{LieAlg},~\mathbf{Mod,}$ $\mathbf{Man}$ and $\mathbf{%
B}^{\mathbf{v}}$ be the category of Lie algebras, modules, manifolds and
vector bundles respectively.

We know that if
\begin{equation*}
\left( E,\pi ,M\right) \in \left\vert \mathbf{B}^{\mathbf{v}}\right\vert ,
\end{equation*}
\begin{equation*}
\Gamma \left( E,\pi ,M\right) =\left\{ u\in \mathbf{Man}\left( M,E\right)
:u\circ \pi =Id_{M}\right\}
\end{equation*}
and
\begin{equation*}
\mathcal{F}\left( M\right) =\mathbf{Man}\left( M,\mathbb{R}\right) ,
\end{equation*}
then $\left( \Gamma \left( E,\pi ,M\right) ,+,\cdot \right) $ is a $\mathcal{%
F}\left( M\right) $-module.

If \ $\left( \varphi ,\varphi _{0}\right) \in \mathbf{B}^{\mathbf{v}}\left(
\left( E,\pi ,M\right) ,\left( E^{\prime },\pi ^{\prime },M^{\prime }\right)
\right) $ such that $\varphi _{0}\in Iso_{\mathbf{Man}}\left( M,M^{\prime
}\right) ,$ then, using the operation
\begin{equation*}
\begin{array}{ccc}
\mathcal{F}\left( M\right) \times \Gamma \left( E^{\prime },\pi ^{\prime
},M^{\prime }\right) & ^{\underrightarrow{~\ \ \cdot ~\ \ }} & \Gamma \left(
E^{\prime },\pi ^{\prime },M^{\prime }\right) \\
\left( f,u^{\prime }\right) & \longmapsto & f\circ \varphi _{0}^{-1}\cdot
u^{\prime }%
\end{array}%
\end{equation*}%
it results that $\left( \Gamma \left( E^{\prime },\pi ^{\prime },M^{\prime
}\right) ,+,\cdot \right) $ is a $\mathcal{F}\left( M\right) $-module and we
obtain the $\mathbf{Mod}$-morphism%
\begin{equation*}
\begin{array}{ccc}
\Gamma \left( E,\pi ,M\right) & ^{\underrightarrow{~\ \ \Gamma \left(
\varphi ,\varphi _{0}\right) ~\ \ }} & \Gamma \left( E^{\prime },\pi
^{\prime },M^{\prime }\right) \\
u & \longmapsto & \Gamma \left( \varphi ,\varphi _{0}\right) u%
\end{array}%
\end{equation*}%
defined by
\begin{equation*}
\begin{array}{c}
\Gamma \left( \varphi ,\varphi _{0}\right) u\left( y\right) =\varphi \left(
u_{\varphi _{0}^{-1}\left( y\right) }\right) =\left( \varphi \circ u\circ
\varphi _{0}^{-1}\right) \left( y\right) ,%
\end{array}%
\end{equation*}%
for any $y\in M^{\prime }.$

If $\left( F,\nu ,M\right) \in \left\vert \mathbf{B}^{\mathbf{v}}\right\vert
$ so that there exists
\begin{equation*}
\begin{array}{c}
\left( \rho ,Id_{M}\right) \in \mathbf{B}^{\mathbf{v}}\left( \left( F,\nu
,M\right) ,\left( TM,\tau _{M},M\right) \right)%
\end{array}%
\end{equation*}%
and an operation
\begin{equation*}
\begin{array}{ccc}
\Gamma \left( F,\nu ,M\right) \times \Gamma \left( F,\nu ,M\right) & ^{%
\underrightarrow{~\ \left[ ,\right] _{F}~\ }} & \Gamma \left( F,\nu ,M\right)
\\
\left( u,v\right) & \longmapsto & \left[ u,v\right] _{F}%
\end{array}%
\end{equation*}%
with the following properties:\bigskip

\noindent $\qquad LA_{1}$. \emph{the equality holds good }%
\begin{equation*}
\begin{array}{c}
\left[ u,f\cdot v\right] _{F}=f\left[ u,v\right] _{F}+\Gamma \left( \rho
,Id_{M}\right) \left( u\right) f\cdot v,%
\end{array}%
\end{equation*}%
\qquad \quad\ \ \emph{for all }$u,v\in \Gamma \left( F,\nu ,M\right) $\emph{%
\ and} $f\in \mathcal{F}\left( M\right) .$

\medskip $LA_{2}$. \emph{the }$4$\emph{-tuple} $\left( \Gamma \left( F,\nu
,M\right) ,+,\cdot ,\left[ ,\right] _{F}\right) $ \emph{is a Lie} $\mathcal{F%
}\left( M\right) $\emph{-algebra,}

$LA_{3}$. \emph{the }$\mathbf{Mod}$\emph{-morphism }$\Gamma \left( \rho
,Id_{M}\right) $\emph{\ is a }$\mathbf{LieAlg}$\emph{-morphism of }%
\begin{equation*}
\left( \Gamma \left( F,\nu ,M\right) ,+,\cdot ,\left[ ,\right] _{F}\right)
\end{equation*}%
\emph{\ source and }%
\begin{equation*}
\left( \Gamma \left( TN,\tau _{N},M\right) ,+,\cdot ,\left[ ,\right]
_{TM}\right)
\end{equation*}%
\emph{target, then the triple }%
\begin{equation*}
\begin{array}{c}
\left( \left( F,\nu ,M\right) ,\left[ ,\right] _{F},\left( \rho
,Id_{M}\right) \right)%
\end{array}%
\leqno(1.1)
\end{equation*}%
\emph{is an object of the category }$\mathbf{LA}$\emph{\ of Lie algebroids}$%
\mathbf{.}$ The couple $\left( \left[ ,\right] _{F},\left( \rho
,Id_{M}\right) \right) $ is called\emph{\ Lie algebroid structure.}

The geometry of the standard Lie algebroid%
\begin{equation*}
\begin{array}{c}
\left( \left( TM,\tau _{M},M\right) ,\left[ ,\right] _{TM},\left(
Id_{TM},Id_{M}\right) \right)
\end{array}%
\leqno(1.2)
\end{equation*}%
where extensively developed with the help of a metrical structure%
\begin{equation*}
\begin{array}{c}
g=g_{ij}\cdot dx^{i}\otimes dx^{j}.%
\end{array}%
\leqno(1.3)
\end{equation*}

A metrical compatible linear connection with free torsion was the
Levi-Civita linear connection. This linear connection, plays an important
role in this geometry, called by use \emph{Riemannian geometry.}

We know that a \emph{smooth Finsler fundamental function} on the tangent
vector bundle $\left( TM,\tau _{M},M\right) $ is a mapping $TM~^{%
\underrightarrow{\ F\ }}~\mathbb{R}_{+}$ which satisfies the following
conditions:\medskip

1. $F\circ u\in C^{\infty }\left( M\right) $, for any $u\in \Gamma \left(
TM,\tau _{M},M\right) \setminus \left\{ 0\right\} $;\smallskip

2. $F\circ 0\in C^{0}\left( M\right) $, where $0$ means the null section of $%
\left( TM,\tau _{M},M\right) $;\smallskip

3. $F$ is positively $1$-homogenous on the fibres of vector bundle $\left(
TM,\tau _{M},M\right) ;$\smallskip

4. For any vector local $m+m$-chart $\left( U,s_{U}\right) $ of $\left(
TM,\tau _{M},M\right) ,$ the hessian:%
\begin{equation*}
\begin{array}{c}
\left\Vert \frac{\partial ^{2}F^{2}\left( u_{x}\right) }{\partial
y^{i}\partial y^{j}}\right\Vert%
\end{array}%
\leqno(1.4)
\end{equation*}%
is positively define for any $u_{x}\in \tau _{M}^{-1}\left( U\right)
\backslash \left\{ 0_{x}\right\} $.

After the Einstein formulation of general relativity, the Riemannian
geometry become widely used and the geometry of the standard Lie algebroid%
\begin{equation*}
\begin{array}{c}
\left( \left( TTM,\tau _{TM},TM\right) ,\left[ ,\right] _{TTM},\left(
Id_{TTM},Id_{TM}\right) \right)%
\end{array}%
\leqno(1.5)
\end{equation*}%
was extensively studied with the help of the metrical structure%
\begin{equation*}
\begin{array}{c}
g=g_{ij}\cdot dx^{i}\otimes dx^{j}+g_{ij}\cdot \delta y^{i}\otimes \delta
y^{j},%
\end{array}%
\leqno(1.6)
\end{equation*}%
where
\begin{equation*}
\begin{array}{c}
g_{ij}=\frac{1}{2}\cdot \frac{\partial ^{2}F^{2}}{\partial y^{i}\partial
y^{j}}%
\end{array}%
\leqno(1.7)
\end{equation*}%
varying smoothly, for any $i,j\in \overline{1,m}.$

This geometry is called the Finsler geometry. Probably, the first work in
Finsler geometry was the PhD thesis of Paul Finsler $(1918).$

First linear connections metrical compatible were proposed by J. L. Synge $%
\left( 1925\right) ,$ J. H. Taylor $\left( 1925\right) ,$ L. Berwald $\left(
1928\right) $ (see: $\left[ 13\right] $). E Cartan $\left( 1934\right) $
(see: $\left[ 18\right] $) propesed a metric compatible linear connection
with the largest number of nonvanishing components of torsion tensor. After
a short time, S. S. Chern $\left[ 19,20\right] $ proposed a different
generalization which is identical with the connection propsed later by Rund
(see $\left[ 5\right] $).

These connections where used to prove many results from Riemannian geometry
in Finslerian geometry context $\left[ 1,11\right] .$ Another useful
connection in Finsler geometry is the Berwald connection $\left[ 11,14,24%
\right] .$ The Berwald connection is torsion-free, but is not metric
compatible. The Berwald curvature tensors are of two types: an $hh-$ one not
unlike the Riemannian curvature tensor and an $hv-$ one which automatically
vanishes in the Riemannian setting. Berwald's connections have been
indispensable to the geometry of path spaces. In $\left[ 2\right] $ and $%
\left[ 28\right] $ one can find caracterizations of these connections
illustrating their similarities and differance.

Entusiasts of metric compatibility where not to be outdone. It is an amusing
irony that although Finsler geometry starts with only a norm in any given
tangent space, it regains an entire family of inner products, one for each
direction in that tangent space. This is why one can still make sense of
metric compatibility in the Finsler setting. The Cartan connection remains
immensely popular with the Matsumoto and Miron schools of Finsler geometry.

Important contributions to the geometry of Finsler spaces were obtained by
M. Abate and G. Patritio $\left[ 1\right] ,$ D. Bao, S. S. Cern and Z. Shen $%
\left[ 11\right] ,$ A. Bejancu $\left[ 12\right] ,$ L. Berwald $\left[ 14%
\right] ,$ H. Busmann $\left[ 16\right] ,$ A. Kawaguchi $\left[ 22\right] ,$
E. Cartan $\left[ 17,18\right] .$

We know that a \emph{smooth Lagrange fundamental function} on the tangent
vector bundle $\left( TM,\tau _{M},M\right) $ is a mapping $E~\ ^{%
\underrightarrow{\ \ L\ \ }}~\ \mathbb{R}$ which satisfies the following
conditions:\medskip

1. $L\circ u\in C^{\infty }\left( M\right) $, for any $u\in \Gamma \left(
E,\pi ,M\right) \setminus \left\{ 0\right\} $;\smallskip

2. $L\circ 0\in C^{0}\left( M\right) $, where $0$ means the null section of $%
\left( E,\pi ,M\right) .$\medskip

If the Hessian matrix with entries%
\begin{equation*}
\begin{array}[b]{c}
g_{ij}=\frac{1}{2}\frac{\partial ^{2}L}{\partial y^{i}\partial y^{j}}%
\end{array}%
\end{equation*}%
is everywhere nondegenerate, then we say that \emph{the Lagrangian is
regular.}

Using a regular Lagrangian, the geometry of the standard Lie algebroid%
\begin{equation*}
\begin{array}{c}
\left( \left( TTM,\tau _{TM},TM\right) ,\left[ ,\right] _{TTM},\left(
Id_{TTM},Id_{TM}\right) \right)%
\end{array}%
\leqno(1.5)^{\prime }
\end{equation*}%
was extensively studied with the help of the metrical structure%
\begin{equation*}
\begin{array}{c}
g=g_{ij}\cdot dx^{i}\otimes dx^{j}+g_{ij}\cdot \delta y^{i}\otimes \delta
y^{j},%
\end{array}%
\leqno(1.6)^{\prime }
\end{equation*}%
where
\begin{equation*}
\begin{array}{c}
g_{ij}=\frac{1}{2}\cdot \frac{\partial ^{2}L}{\partial y^{i}\partial y^{j}}%
\end{array}%
\leqno(1.7)^{\prime }
\end{equation*}%
varying smoothly, for any $i,j\in \overline{1,m}$ and for any local vector $%
m+m$-chart $\left( U,s_{U}\right) $ the matrix $\left\Vert g_{ij}\right\Vert
$ has constant signature on $\tau _{M}^{-1}\left( U\right) -\left\{
0_{x}\right\} .$

An important generalization of Finsler geometry, called by use the\emph{\
Lagrange geometry,} have been proposed. This geometry can be developed by
the methods of Finsler geometry, but using the main mechanical properties of
the Lagrangians. The notion of Lagrange space were introduced and studied by
J. Kern $\left[ 23\right] $ and R. Miron $\left[ 25,27\right] .$ The
nonlinear connections and the distinguished linear connections which depend
only on Lagrangian was presented in the framework of Lagrange spaces. The
geometry of Lagrange spaces have been developed in many proceedings and
monographs $\left[ 10,11,15,24,25,28,29,30,31,32\right] .$The Roumanian
school initiated by R. Miron has important contributions $\left[ 25,29,33%
\right] $.

A natural generalization of Lagrange geometry is provided by the notion of
\emph{generalized Lagrange space }introduced by R. Miron in $\left[ 26\right]
$. In this general framework, the geometry of the standard Lie algebroid%
\begin{equation*}
\begin{array}{c}
\left( \left( TTM,\tau _{TM},TM\right) ,\left[ ,\right] _{TTM},\left(
Id_{TTM},Id_{TM}\right) \right)%
\end{array}%
\leqno(1.5)^{\prime \prime }
\end{equation*}%
continous to be studied with the help of a metrical structure%
\begin{equation*}
\begin{array}{c}
g=g_{ij}\cdot dx^{i}\otimes dx^{j}+g_{ij}\cdot \delta y^{i}\otimes \delta
y^{j},%
\end{array}%
\leqno(1.6)^{\prime \prime }
\end{equation*}%
where
\begin{equation*}
\begin{array}{c}
g_{ij}\cdot dy^{i}\otimes dy^{j}\in \mathcal{T}~_{2}^{0}\left( VTTM,{\tau
_{TM}},TM\right)%
\end{array}%
\leqno(1.7)^{\prime \prime }
\end{equation*}

The generalized Lagrange geometry constructed by R. Miron starts with the
problem of associating of a cononical connection to these spaces. If $\Gamma
$ is a connection for the tangent vector bundle $\left( TM,\tau
_{M},M\right) ,$ then one determines a metrical distinguished linear
connection uniquely in some conditions. Its coefficients are expresed by the
generalized Christoffel symbols. Important results are presented by M.
Anastasiei in $\left[ 3,4,6\right] .$

Using the notion of generalized Lie algebroid $\left[ 8\right] $, in the
paper $\left[ 7\right] ,$ we started the study of the geometry of the Lie
algebroid generalized tangent bundle
\begin{equation*}
\begin{array}{c}
\left( \left( \left( \rho ,\eta \right) TE,\left( \rho ,\eta \right) \tau
_{E},E\right) ,\left[ ,\right] _{\left( \rho ,\eta \right) TE},\left( \tilde{%
\rho},Id_{E}\right) \right)%
\end{array}%
\end{equation*}%
of a vector bundle $\left( E,\pi ,M\right) ,$ from mechanical $\left( \rho
,\eta \right) $-systems point of view. A complet theory of $\left( \rho
,\eta \right) $-(semi)sprays and a Lagrangian formalism for Lagrange
mechanical $\left( \rho ,\eta \right) $-systems are presented. In
particular, if the diffeomorphisms used are identities, then the classical
results are obtained.

In this paper we continous to study the geometry of the Lie algebroid
generalized tangent bundle from the metrizability point of view. Using the
generalized Lie algebroid, in Section $2,$ we build the Lie algebroid
generalized tangent bundle. Using the vertical interior differential system
(see $\left[ 9\right] $), the $\left( \rho ,\eta \right) $-connections
theory and adapted basis are presenten in Section $3$. The distinguished
linear $\left( \rho ,\eta \right) $-connections are presented in Section $4.$
In Section $5,$ we develop the theory of metrizability of the Lie algebroid
generalized tangent bundle and a lot of important results are presented.
Notice that our theory is a progress, because in particular, if the
diffeomorphisms used are identities, then we obtain important results in the
framework of Lie algebroids. In the particular case of the Lie algebroid
tangent bundle
\begin{equation*}
\begin{array}{c}
\left( \left( TE,\tau _{E},E\right) ,\left[ ,\right] _{TE},\left(
Id_{TE},Id_{E}\right) \right)%
\end{array}%
\end{equation*}%
of a vector bundle $\left( E,\pi ,M\right) ,$ we obtain all classical
results. (see $\left[ 29,31,33\right] $).

Finally, in Section $6,$ a new class of Lagrange spaces, called by use,
generalized Lagrange $\left( \rho ,\eta \right) $-spaces, Lagrange $\left(
\rho ,\eta \right) $-spaces and Finsler $\left( \rho ,\eta \right) $-spaces
are presented. In the particular case of Lie algebroids, new and important
results are obtained. When all morphisms are identities, then the classical
results for generalized Lagrange spaces are obtained.

As the Lagrange $\left( Id_{TM},Id_{M}\right) $-spaces and Finsler $\left(
Id_{TM},Id_{M}\right) $-spaces are the usual Lagrange spaces and Finsler
spaces, then we ask:

\begin{quote}
{\small - Can we study the geometry of the Lie algebroid generalized tangent
bundle }%
\begin{equation*}
\begin{array}{c}
\left( \left( \left( \rho ,\eta \right) TE,\left( \rho ,\eta \right) \tau
_{E},E\right) ,\left[ ,\right] _{\left( \rho ,\eta \right) TE},\left( \tilde{%
\rho},Id_{E}\right) \right)%
\end{array}%
\end{equation*}%
{\small from the Lagrange }$\left( \rho ,\eta \right) ${\small -spaces and
Finsler }$\left( \rho ,\eta \right) ${\small -spaces point of view? In
particular, if all morphisms are identities, then are obtain the classical
results for Lagrange and Finsler? }
\end{quote}

This is a new direction by reserch of the geometry of the Lie algebroid
generalized tangent bundle.

\section{The Lie algebroid generalized tangent bundle}

\ \ \ \

If $M,N\in \left\vert \mathbf{Man}\right\vert ,$ $h\in Iso_{\mathbf{Man}%
}\left( M,N\right) $, $\eta \in Iso_{\mathbf{Man}}\left( N,M\right) $ and $%
\left( F,\nu ,N\right) \in \left\vert \mathbf{B}^{\mathbf{v}}\right\vert $
so that there exists
\begin{equation*}
\begin{array}{c}
\left( \rho ,\eta \right) \in \mathbf{B}^{\mathbf{v}}\left( \left( F,\nu
,N\right) ,\left( TM,\tau _{M},M\right) \right)%
\end{array}%
\end{equation*}%
and an operation
\begin{equation*}
\begin{array}{ccc}
\Gamma \left( F,\nu ,N\right) \times \Gamma \left( F,\nu ,N\right) & ^{%
\underrightarrow{\left[ ,\right] _{F,h}}} & \Gamma \left( F,\nu ,N\right) \\
\left( u,v\right) & \longmapsto & \left[ u,v\right] _{F,h}%
\end{array}%
\end{equation*}%
with the following properties:\bigskip

\noindent $\qquad GLA_{1}$. \emph{the equality holds good }%
\begin{equation*}
\begin{array}{c}
\left[ u,f\cdot v\right] _{F,h}=f\left[ u,v\right] _{F,h}+\Gamma \left(
Th\circ \rho ,h\circ \eta \right) \left( u\right) f\cdot v,%
\end{array}%
\end{equation*}%
\qquad \quad\ \ \emph{for all }$u,v\in \Gamma \left( F,\nu ,N\right) $\emph{%
\ and} $f\in \mathcal{F}\left( N\right) .$

\medskip $GLA_{2}$. \emph{the }$4$\emph{-tuple} $\left( \Gamma \left( F,\nu
,N\right) ,+,\cdot ,\left[ ,\right] _{F,h}\right) $ \emph{is a Lie} $%
\mathcal{F}\left( N\right) $\emph{-algebra,}

$GLA_{3}$. \emph{the }$\mathbf{Mod}$\emph{-morphism }$\Gamma \left( Th\circ
\rho ,h\circ \eta \right) $\emph{\ is a }$\mathbf{LieAlg}$\emph{-morphism of
}%
\begin{equation*}
\left( \Gamma \left( F,\nu ,N\right) ,+,\cdot ,\left[ ,\right] _{F,h}\right)
\end{equation*}%
\emph{\ source and }%
\begin{equation*}
\left( \Gamma \left( TN,\tau _{N},N\right) ,+,\cdot ,\left[ ,\right]
_{TN}\right)
\end{equation*}%
\emph{target, then the triple }%
\begin{equation*}
\begin{array}{c}
\left( \left( F,\nu ,N\right) ,\left[ ,\right] _{F,h},\left( \rho ,\eta
\right) \right)%
\end{array}%
\leqno(2.1)
\end{equation*}%
\emph{is an object of the category }$\mathbf{GLA}$\emph{\ of generalized Lie
algebroids}$\mathbf{.}$ The couple $\left( \left[ ,\right] _{F,h},\left(
\rho ,\eta \right) \right) $ is called \emph{generalized Lie algebroid
structure. }(see $\left[ 8\right] $)

In particular, if $\left( \eta ,h\right) =\left( Id_{M},Id_{M}\right) ,$\ we
obtain the definition of the Lie algebroid.

\begin{itemize}
\item Locally, for any $\alpha ,\beta \in \overline{1,p},$ we set $\left[
t_{\alpha },t_{\beta }\right] _{F,h}=L_{\alpha \beta }^{\gamma }t_{\gamma }.$
We easily obtain that $L_{\alpha \beta }^{\gamma }=-L_{\beta \alpha
}^{\gamma },~$for any $\alpha ,\beta ,\gamma \in \overline{1,p}.$
\end{itemize}

The real local functions $L_{\alpha \beta }^{\gamma },~\alpha ,\beta ,\gamma
\in \overline{1,p}$ are called the \emph{structure functions of the
generalized Lie algebroid }$\left( \left( F,\nu ,N\right) ,\left[ ,\right]
_{F,h},\left( \rho ,\eta \right) \right) .$

\begin{itemize}
\item We assume the following diagrams:%
\begin{equation*}
\begin{array}[b]{ccccc}
F & ^{\underrightarrow{~\ \ \ \rho ~\ \ }} & TM & ^{\underrightarrow{~\ \ \
Th~\ \ }} & TN \\
~\downarrow \nu &  & ~\ \ \ \downarrow \tau _{M} &  & ~\ \ \ \downarrow \tau
_{N} \\
N & ^{\underrightarrow{~\ \ \ \eta ~\ \ }} & M & ^{\underrightarrow{~\ \ \
h~\ \ }} & N \\
&  &  &  &  \\
\left( \chi ^{\tilde{\imath}},z^{\alpha }\right) &  & \left(
x^{i},y^{i}\right) &  & \left( \chi ^{\tilde{\imath}},z^{\tilde{\imath}%
}\right)%
\end{array}%
\end{equation*}

where $i,\tilde{\imath}\in \overline{1,m}$ and $\alpha \in \overline{1,p}.$
\end{itemize}

If%
\begin{equation*}
\left( \chi ^{\tilde{\imath}},z^{\alpha }\right) \longrightarrow \left( \chi
^{\tilde{\imath}\prime }\left( \chi ^{\tilde{\imath}}\right) ,z^{\alpha
\prime }\left( \chi ^{\tilde{\imath}},z^{\alpha }\right) \right) ,
\end{equation*}%
\begin{equation*}
\left( x^{i},y^{i}\right) \longrightarrow \left( x^{i%
%TCIMACRO{\U{b4}}%
%BeginExpansion
{\acute{}}%
%EndExpansion
}\left( x^{i}\right) ,y^{i%
%TCIMACRO{\U{b4}}%
%BeginExpansion
{\acute{}}%
%EndExpansion
}\left( x^{i},y^{i}\right) \right)
\end{equation*}%
and
\begin{equation*}
\left( \chi ^{\tilde{\imath}},z^{\tilde{\imath}}\right) \longrightarrow
\left( \chi ^{\tilde{\imath}\prime }\left( \chi ^{\tilde{\imath}}\right) ,z^{%
\tilde{\imath}\prime }\left( \chi ^{\tilde{\imath}},z^{\tilde{\imath}%
}\right) \right) ,
\end{equation*}%
then
\begin{equation*}
\begin{array}[b]{c}
z^{\alpha
%TCIMACRO{\U{b4}}%
%BeginExpansion
{\acute{}}%
%EndExpansion
}=\Lambda _{\alpha }^{\alpha
%TCIMACRO{\U{b4}}%
%BeginExpansion
{\acute{}}%
%EndExpansion
}z^{\alpha }%
\end{array}%
,
\end{equation*}%
\begin{equation*}
\begin{array}[b]{c}
y^{i%
%TCIMACRO{\U{b4}}%
%BeginExpansion
{\acute{}}%
%EndExpansion
}=\frac{\partial x^{i%
%TCIMACRO{\U{b4}}%
%BeginExpansion
{\acute{}}%
%EndExpansion
}}{\partial x^{i}}y^{i}%
\end{array}%
\end{equation*}%
and
\begin{equation*}
\begin{array}{c}
z^{\tilde{\imath}\prime }=\frac{\partial \chi ^{\tilde{\imath}\prime }}{%
\partial \chi ^{\tilde{\imath}}}z^{\tilde{\imath}}.%
\end{array}%
\end{equation*}

We assume that $\left( \theta ,\mu \right) \overset{put}{=}\left( Th\circ
\rho ,h\circ \eta \right) $. If $z^{\alpha }t_{\alpha }\in \Gamma \left(
F,\nu ,N\right) $ is arbitrary, then
\begin{equation*}
\begin{array}[t]{l}
\displaystyle%
\begin{array}{c}
\Gamma \left( Th\circ \rho ,h\circ \eta \right) \left( z^{\alpha }t_{\alpha
}\right) f\left( h\circ \eta \left( \varkappa \right) \right) =\vspace*{1mm}
\\
=\left( \theta _{\alpha }^{\tilde{\imath}}z^{\alpha }\frac{\partial f}{%
\partial \varkappa ^{\tilde{\imath}}}\right) \left( h\circ \eta \left(
\varkappa \right) \right) =\left( \left( \rho _{\alpha }^{i}\circ h\right)
\left( z^{\alpha }\circ h\right) \frac{\partial f\circ h}{\partial x^{i}}%
\right) \left( \eta \left( \varkappa \right) \right) ,%
\end{array}%
\end{array}%
\leqno(2.2)
\end{equation*}%
for any $f\in \mathcal{F}\left( N\right) $ and $\varkappa \in N.$

The coefficients $\rho _{\alpha }^{i}$ respectively $\theta _{\alpha }^{%
\tilde{\imath}}$ change to $\rho _{\alpha
%TCIMACRO{\U{b4}}%
%BeginExpansion
{\acute{}}%
%EndExpansion
}^{i%
%TCIMACRO{\U{b4}}%
%BeginExpansion
{\acute{}}%
%EndExpansion
}$ respectively $\theta _{\alpha
%TCIMACRO{\U{b4}}%
%BeginExpansion
{\acute{}}%
%EndExpansion
}^{\tilde{\imath}%
%TCIMACRO{\U{b4}}%
%BeginExpansion
{\acute{}}%
%EndExpansion
}$ according to the rule:
\begin{equation*}
\begin{array}{c}
\rho _{\alpha
%TCIMACRO{\U{b4}}%
%BeginExpansion
{\acute{}}%
%EndExpansion
}^{i%
%TCIMACRO{\U{b4}}%
%BeginExpansion
{\acute{}}%
%EndExpansion
}=\Lambda _{\alpha
%TCIMACRO{\U{b4}}%
%BeginExpansion
{\acute{}}%
%EndExpansion
}^{\alpha }\rho _{\alpha }^{i}\displaystyle\frac{\partial x^{i%
%TCIMACRO{\U{b4}}%
%BeginExpansion
{\acute{}}%
%EndExpansion
}}{\partial x^{i}},%
\end{array}%
\leqno(2.3)
\end{equation*}%
respectively%
\begin{equation*}
\begin{array}{c}
\theta _{\alpha
%TCIMACRO{\U{b4}}%
%BeginExpansion
{\acute{}}%
%EndExpansion
}^{\tilde{\imath}%
%TCIMACRO{\U{b4}}%
%BeginExpansion
{\acute{}}%
%EndExpansion
}=\Lambda _{\alpha
%TCIMACRO{\U{b4}}%
%BeginExpansion
{\acute{}}%
%EndExpansion
}^{\alpha }\theta _{\alpha }^{\tilde{\imath}}\displaystyle\frac{\partial
\varkappa ^{\tilde{\imath}%
%TCIMACRO{\U{b4}}%
%BeginExpansion
{\acute{}}%
%EndExpansion
}}{\partial \varkappa ^{\tilde{\imath}}},%
\end{array}%
\leqno(2.4)
\end{equation*}%
where
\begin{equation*}
\left\Vert \Lambda _{\alpha
%TCIMACRO{\U{b4}}%
%BeginExpansion
{\acute{}}%
%EndExpansion
}^{\alpha }\right\Vert =\left\Vert \Lambda _{\alpha }^{\alpha
%TCIMACRO{\U{b4}}%
%BeginExpansion
{\acute{}}%
%EndExpansion
}\right\Vert ^{-1}.
\end{equation*}

\emph{Remark 2.1 }\emph{The following equalities hold good:}%
\begin{equation*}
\begin{array}{c}
\displaystyle\rho _{\alpha }^{i}\circ h\frac{\partial f\circ h}{\partial
x^{i}}=\left( \theta _{\alpha }^{\tilde{\imath}}\frac{\partial f}{\partial
\varkappa ^{\tilde{\imath}}}\right) \circ h,\forall f\in \mathcal{F}\left(
N\right) .%
\end{array}%
\leqno(2.5)
\end{equation*}%
\emph{and }%
\begin{equation*}
\begin{array}{c}
\displaystyle\left( L_{\alpha \beta }^{\gamma }\circ h\right) \left( \rho
_{\gamma }^{k}\circ h\right) =\left( \rho _{\alpha }^{i}\circ h\right) \frac{%
\partial \left( \rho _{\beta }^{k}\circ h\right) }{\partial x^{i}}-\left(
\rho _{\beta }^{j}\circ h\right) \frac{\partial \left( \rho _{\alpha
}^{k}\circ h\right) }{\partial x^{j}}.%
\end{array}%
\leqno(2.6)
\end{equation*}

\emph{In the particular case of Lie algebroids, the relation }$\left(
2.6\right) $\emph{\ becomes}
\begin{equation*}
\begin{array}{c}
\displaystyle L_{\alpha \beta }^{\gamma }\rho _{\gamma }^{k}=\rho _{\alpha
}^{i}\frac{\partial \rho _{\beta }^{k}}{\partial x^{i}}-\rho _{\beta }^{j}%
\frac{\partial \rho _{\alpha }^{k}}{\partial x^{j}}.%
\end{array}%
\leqno(2.6^{\prime })
\end{equation*}

\textbf{Example 2.1} \emph{Let} $M,N\in \left\vert \mathbf{Man}\right\vert ,$
$h\in Iso_{\mathbf{Man}}\left( M,N\right) $ \emph{and} $\eta \in Iso_{%
\mathbf{Man}}\left( N,M\right) $\emph{\ be. Using the tangent }$\mathbf{B}^{%
\mathbf{v}}$\emph{-morphism }$\left( T\eta ,\eta \right) $\emph{\ and the
operation }%
\begin{equation*}
\begin{array}{ccc}
\Gamma \left( TN,\tau _{N},N\right) \times \Gamma \left( TN,\tau
_{N},N\right) & ^{\underrightarrow{~\ \ \left[ ,\right] _{TN,h}~\ \ }} &
\Gamma \left( TN,\tau _{N},N\right) \\
\left( u,v\right) & \longmapsto & \ \left[ u,v\right] _{TN,h}%
\end{array}%
\end{equation*}%
\emph{where }%
\begin{equation*}
\left[ u,v\right] _{TN,h}=\Gamma \left( T\left( h\circ \eta \right)
^{-1},\left( h\circ \eta \right) ^{-1}\right) \left( \left[ \Gamma \left(
T\left( h\circ \eta \right) ,h\circ \eta \right) u,\Gamma \left( T\left(
h\circ \eta \right) ,h\circ \eta \right) v\right] _{TN}\right) ,
\end{equation*}%
\emph{for any }$u,v\in \Gamma \left( TN,\tau _{N},N\right) $\emph{, we
obtain that}%
\begin{equation*}
\begin{array}{c}
\left( \left( TN,\tau _{N},N\right) ,\left( T\eta ,\eta \right) ,\left[ ,%
\right] _{TN,h}\right)%
\end{array}%
\end{equation*}%
\emph{is a generalized Lie algebroid.}

For any $\mathbf{Man}$-isomorphisms $\eta $ and $h,$ new and interesting
generalized Lie algebroid structures for the tangent vector bundle $\left(
TN,\tau _{N},N\right) $ are obtained$.$ For any base $\left\{ t_{\alpha
},~\alpha \in \overline{1,m}\right\} $ of the module of sections $\left(
\Gamma \left( TN,\tau _{N},N\right) ,+,\cdot \right) $ we obtain the
structure functions%
\begin{equation*}
\begin{array}{c}
L_{\alpha \beta }^{\gamma }=\left( \theta _{\alpha }^{i}\frac{\partial
\theta _{\beta }^{j}}{\partial x^{i}}-\theta _{\beta }^{i}\frac{\partial
\theta _{\alpha }^{j}}{\partial x^{i}}\right) \tilde{\theta}_{j}^{\gamma
},~\alpha ,\beta ,\gamma \in \overline{1,m}%
\end{array}%
\end{equation*}%
where
\begin{equation*}
\theta _{\alpha }^{i},~i,\alpha \in \overline{1,m}
\end{equation*}%
are real local functions so that
\begin{equation*}
\begin{array}{c}
\Gamma \left( T\left( h\circ \eta \right) ,h\circ \eta \right) \left(
t_{\alpha }\right) =\theta _{\alpha }^{i}\frac{\partial }{\partial x^{i}}%
\end{array}%
\end{equation*}%
and
\begin{equation*}
\tilde{\theta}_{j}^{\gamma },~i,\gamma \in \overline{1,m}
\end{equation*}%
are real local functions so that
\begin{equation*}
\begin{array}{c}
\Gamma \left( T\left( h\circ \eta \right) ^{-1},\left( h\circ \eta \right)
^{-1}\right) \left( \frac{\partial }{\partial x^{j}}\right) =\tilde{\theta}%
_{j}^{\gamma }t_{\gamma }.%
\end{array}%
\end{equation*}

In particular, using arbitrary isometries (symmetries, translations,
rotations,...) for the Euclidean $3$-dimensional space $\Sigma ,$ and
arbitrary basis for the module of sections we obtain a lot of generalized
Lie algebroid structures for the tangent vector bundle $\left( T\Sigma ,\tau
_{\Sigma },\Sigma \right) $.

Let $\left( E,\pi ,M\right) $ be a vector bundle. We obtain the $\mathbf{B}^{%
\mathbf{v}}$-morphism%
\begin{equation*}
\begin{array}{ccc}
~\ \ \ \ \ \ \ \ \ \ \ \ \ \pi ^{\ast }\ \left( h^{\ast }F\right) &
\hookrightarrow & F \\
\pi ^{\ast }\left( h^{\ast }\nu \right) \downarrow &  & ~\downarrow \nu \\
~\ \ \ \ \ \ \ \ \ \ \ \ \ \ E & ^{\underrightarrow{~\ \ h\circ \pi ~\ \ }}
& M%
\end{array}%
\leqno(2.7)
\end{equation*}

We take $\left( x^{i},y^{a}\right) $ as canonical local coordinates on $%
\left( E,\pi ,M\right) ,$ where $i\in \overline{1,m}$ and $a\in \overline{1,r%
}.$ Let
\begin{equation*}
\left( x^{i},y^{a}\right) \longrightarrow \left( x^{i%
%TCIMACRO{\U{b4}}%
%BeginExpansion
{\acute{}}%
%EndExpansion
}\left( x^{i}\right) ,y^{a%
%TCIMACRO{\U{b4}}%
%BeginExpansion
{\acute{}}%
%EndExpansion
}\left( x^{i},y^{a}\right) \right)
\end{equation*}%
be a change of coordinates on $\left( E,\pi ,M\right) $. Then the
coordinates $y^{a}$ change to $y^{a%
%TCIMACRO{\U{b4}}%
%BeginExpansion
{\acute{}}%
%EndExpansion
}$ according to the rule:
\begin{equation*}
\begin{array}{c}
y^{a%
%TCIMACRO{\U{b4}}%
%BeginExpansion
{\acute{}}%
%EndExpansion
}=\displaystyle M_{a}^{a%
%TCIMACRO{\U{b4}}%
%BeginExpansion
{\acute{}}%
%EndExpansion
}y^{a}.%
\end{array}%
\leqno(2.8)
\end{equation*}

\textbf{Theorem 2.1 }\emph{Let} $\Big({\overset{\pi ^{\ast }\ \left( h^{\ast
}F\right) }{\rho }},Id_{E}\Big)$ \emph{be the }$\mathbf{B}^{\mathbf{v}}$%
\emph{-morphism of }$\left( \pi ^{\ast }\ \left( h^{\ast }F\right) ,\pi
^{\ast }\left( h^{\ast }\nu \right) ,E\right) $\ \emph{source and} $\left(
TE,\tau _{E},E\right) $\ \emph{target, where}%
\begin{equation*}
\begin{array}{rcl}
\ \pi ^{\ast }\ \left( h^{\ast }F\right) & ^{\underrightarrow{\overset{\pi
^{\ast }\ \left( h^{\ast }F\right) }{\rho }}} & TE \\
\displaystyle Z^{\alpha }T_{\alpha }\left( u_{x}\right) & \longmapsto & %
\displaystyle\left( Z^{\alpha }\cdot \rho _{\alpha }^{i}\circ h\circ \pi
\right) \frac{\partial }{\partial x^{i}}\left( u_{x}\right)%
\end{array}%
\leqno(2.9)
\end{equation*}

\emph{Using the operation}
\begin{equation*}
\begin{array}{ccc}
\Gamma \left( \pi ^{\ast }\ \left( h^{\ast }F\right) ,\pi ^{\ast }\left(
h^{\ast }\nu \right) ,E\right) ^{2} & ^{\underrightarrow{~\ \ \left[ ,\right]
_{\pi ^{\ast }\ \left( h^{\ast }F\right) }~\ \ }} & \Gamma \left( \pi ^{\ast
}\ \left( h^{\ast }F\right) ,\pi ^{\ast }\left( h^{\ast }\nu \right)
,E\right)%
\end{array}%
\end{equation*}%
\emph{defined by}%
\begin{equation*}
\begin{array}{ll}
\left[ T_{\alpha },T_{\beta }\right] _{\pi ^{\ast }\ \left( h^{\ast
}F\right) } & =L_{\alpha \beta }^{\gamma }\circ h\circ \pi \cdot T_{\gamma },%
\vspace*{1mm} \\
\left[ T_{\alpha },fT_{\beta }\right] _{\pi ^{\ast }\ \left( h^{\ast
}F\right) } & \displaystyle=fL_{\alpha \beta }^{\gamma }\circ h\circ \pi
T_{\gamma }+\rho _{\alpha }^{i}\circ h\circ \pi \frac{\partial f}{\partial
x^{i}}T_{\beta },\vspace*{1mm} \\
\left[ fT_{\alpha },T_{\beta }\right] _{\pi ^{\ast }\ \left( h^{\ast
}F\right) } & =-\left[ T_{\beta },fT_{\alpha }\right] _{\pi ^{\ast }\ \left(
h^{\ast }F\right) },%
\end{array}%
\leqno(2.10)
\end{equation*}%
\emph{for any} $f\in \mathcal{F}\left( E\right) ,$ \emph{it results that}
\begin{equation*}
\begin{array}{c}
\left( \left( \pi ^{\ast }\ \left( h^{\ast }F\right) ,\pi ^{\ast }\left(
h^{\ast }\nu \right) ,E\right) ,\left[ ,\right] _{\pi ^{\ast }\ \left(
h^{\ast }F\right) },\left( \overset{\pi ^{\ast }\ \left( h^{\ast }F\right) }{%
\rho },Id_{E}\right) \right)%
\end{array}%
\end{equation*}%
\emph{is a Lie algebroid which is called the pull-back Lie algebroid of the
generalized Lie algebroid }$\left( \left( F,\nu ,N\right) ,\left[ ,\right]
_{F,h},\left( \rho ,\eta \right) \right) $

If $z=z^{\alpha }t_{\alpha }\in \Gamma \left( F,\nu ,N\right) ,$ then we
obtain the section%
\begin{equation*}
Z=\left( z^{\alpha }\circ h\circ \pi \right) T_{\alpha }\in \Gamma \left(
\pi ^{\ast }\left( h^{\ast }F\right) ,\pi ^{\ast }\left( h^{\ast }\nu
\right) ,E\right)
\end{equation*}%
so that $Z\left( u_{x}\right) =z\left( h\left( x\right) \right) ,$ for any $%
u_{x}\in \pi ^{-1}\left( U{\cap h}^{-1}V\right) .$

Let
\begin{equation*}
\begin{array}[b]{c}
\left( \partial _{i},\dot{\partial}_{a}\right) \overset{put}{=}\left( \frac{%
\partial }{\partial x^{i}},\frac{\partial }{\partial y^{a}}\right)%
\end{array}%
\end{equation*}
be the base sections for the Lie $\mathcal{F}\left( E\right) $-algebra
\begin{equation*}
\left( \Gamma \left( TE,\tau _{E},E\right) ,+,\cdot ,\left[ ,\right]
_{TE}\right) .
\end{equation*}

For any sections%
\begin{equation*}
\begin{array}{c}
Z^{\alpha }T_{\alpha }\in \Gamma \left( \pi ^{\ast }\left( h^{\ast }F\right)
,\pi ^{\ast }\left( h^{\ast }F\right) ,E\right)%
\end{array}%
\end{equation*}%
and%
\begin{equation*}
\begin{array}{c}
Y^{a}\dot{\partial}_{a}\in \Gamma \left( VTE,\tau _{E},E\right)%
\end{array}%
\end{equation*}%
we obtain the section
\begin{equation*}
\begin{array}{c}
Z^{\alpha }\tilde{\partial}_{\alpha }+Y^{a}\overset{\cdot }{\tilde{\partial}}%
_{a}=:Z^{\alpha }\left( T_{\alpha }\oplus \left( \rho _{\alpha }^{i}\circ
h\circ \pi \right) \partial _{i}\right) +Y^{a}\left( 0_{\pi ^{\ast }\left(
h^{\ast }F\right) }\oplus \dot{\partial}_{a}\right) \vspace*{1mm} \\
=Z^{\alpha }T_{\alpha }\oplus \left( Z^{\alpha }\left( \rho _{\alpha
}^{i}\circ h\circ \pi \right) \partial _{i}+Y^{a}\dot{\partial}_{a}\right)
\in \Gamma \left( \pi ^{\ast }\left( h^{\ast }F\right) \oplus TE,\overset{%
\oplus }{\pi },E\right) .%
\end{array}%
\end{equation*}

Since we have
\begin{equation*}
\begin{array}{c}
Z^{\alpha }\displaystyle\tilde{\partial}_{\alpha }+Y^{a}\overset{\cdot }{%
\tilde{\partial}}_{a}=0 \\
\Updownarrow \\
Z^{\alpha }T_{\alpha }=0~\wedge Z^{\alpha }\left( \rho _{\alpha }^{i}\circ
h\circ \pi \right) \partial _{i}+Y^{a}\dot{\partial}_{a}=0,%
\end{array}%
\end{equation*}%
it implies $Z^{\alpha }=0,~\alpha \in \overline{1,p}$ and $Y^{a}=0,~a\in
\overline{1,r}.$

Therefore, the sections $\tilde{\partial}_{1},...,\tilde{\partial}_{p},%
\overset{\cdot }{\tilde{\partial}}_{1},...,\overset{\cdot }{\tilde{\partial}}%
_{r}$ are linearly independent.\smallskip

We consider the vector subbundle $\left( \left( \rho ,\eta \right) TE,\left(
\rho ,\eta \right) \tau _{E},E\right) $ of the vector bundle$\left( \pi
^{\ast }\left( h^{\ast }F\right) \oplus TE,\overset{\oplus }{\pi },E\right)
, $ for which the $\mathcal{F}\left( E\right) $-module of sections is the $%
\mathcal{F}\left( E\right) $-submodule of $\left( \Gamma \left( \pi ^{\ast
}\left( h^{\ast }F\right) \oplus TE,\overset{\oplus }{\pi },E\right)
,+,\cdot \right) ,$ generated by the set of sections $\left( \tilde{\partial}%
_{\alpha },\overset{\cdot }{\tilde{\partial}}_{a}\right) .$

The base sections $\left( \tilde{\partial}_{\alpha },\overset{\cdot }{\tilde{%
\partial}}_{a}\right) $ will be called the \emph{natural }$\left( \rho ,\eta
\right) $\emph{-base.}

The matrix of coordinate transformation on $\left( \left( \rho ,\eta \right)
TE,\left( \rho ,\eta \right) \tau _{E},E\right) $ at a change of fibred
charts is
\begin{equation*}
\left\Vert
\begin{array}{cc}
\Lambda _{\alpha }^{\alpha
%TCIMACRO{\U{b4}}%
%BeginExpansion
{\acute{}}%
%EndExpansion
}\circ h\circ \pi & 0\vspace*{1mm} \\
\left( \rho _{a}^{i}\circ h\circ \pi \right) \displaystyle\frac{\partial
M_{b}^{a%
%TCIMACRO{\U{b4}}%
%BeginExpansion
{\acute{}}%
%EndExpansion
}\circ \pi }{\partial x_{i}}y^{b} & M_{a}^{a%
%TCIMACRO{\U{b4}}%
%BeginExpansion
{\acute{}}%
%EndExpansion
}\circ \pi%
\end{array}%
\right\Vert .\leqno(2.11)
\end{equation*}

Easily we obtain

\textbf{Theorem 2.2 }\emph{Let} $\left( \tilde{\rho},Id_{E}\right) $\ \emph{%
be the} $\mathbf{B}^{\mathbf{v}}$\emph{-morphism of }$\left( \left( \rho
,\eta \right) TE,\left( \rho ,\eta \right) \tau _{E},E\right) $\ \emph{%
source and }$\left( TE,\tau _{E},E\right) $\ \emph{target, where}
\begin{equation*}
\begin{array}{rcl}
\left( \rho ,\eta \right) TE\!\!\! & \!\!^{\underrightarrow{\tilde{\ \ \rho
\ \ }}}\!\!\! & \!\!TE\vspace*{2mm} \\
\left( Z^{\alpha }\tilde{\partial}_{\alpha }+Y^{a}\overset{\cdot }{\tilde{%
\partial}}_{a}\right) \!(u_{x})\!\!\!\! & \!\!\longmapsto \!\!\! &
\!\!\left( \!Z^{\alpha }\!\left( \rho _{\alpha }^{i}{\circ }h{\circ }\pi
\!\right) \!\partial _{i}{+}Y^{a}\dot{\partial}_{a}\right) \!(u_{x})\!\!.%
\end{array}%
\leqno(2.12)
\end{equation*}

\emph{Using the operation}
\begin{equation*}
\begin{array}{ccc}
\Gamma \left( \left( \rho ,\eta \right) TE,\left( \rho ,\eta \right) \tau
_{E},E\right) ^{2} & ^{\underrightarrow{~\ \ \left[ ,\right] _{\left( \rho
,\eta \right) TE}~\ \ }} & \Gamma \left( \left( \rho ,\eta \right) TE,\left(
\rho ,\eta \right) \tau _{E},E\right)%
\end{array}%
\end{equation*}%
\emph{defined by}%
\begin{equation*}
\begin{array}{l}
\left[ Z_{1}^{\alpha }\tilde{\partial}_{\alpha }+Y_{1}^{a}\overset{\cdot }{%
\tilde{\partial}}_{a},Z_{2}^{\beta }\tilde{\partial}_{\beta }+Y_{2}^{b}%
\overset{\cdot }{\tilde{\partial}}_{b}\right] _{\left( \rho ,\eta \right) TE}%
\vspace*{1mm} \\
\displaystyle=\left[ Z_{1}^{\alpha }T_{\alpha },Z_{2}^{\beta }T_{\beta }%
\right] _{\pi ^{\ast }\left( h^{\ast }F\right) }\oplus \left[ Z_{1}^{\alpha
}\left( \rho _{\alpha }^{i}\circ h\circ \pi \right) \partial _{i}+Y_{1}^{a}%
\dot{\partial}_{a},\right. \vspace*{1mm} \\
\hfill \displaystyle\left. Z_{2}^{\beta }\left( \rho _{\beta }^{j}\circ
h\circ \pi \right) \partial _{j}+Y_{2}^{b}\dot{\partial}_{b}\right] _{TE},%
\end{array}%
\leqno(2.13)
\end{equation*}%
\emph{for any} $Z_{1}^{\alpha }\tilde{\partial}_{\alpha }+Y_{1}^{a}\overset{%
\cdot }{\tilde{\partial}}_{a}$\emph{\ and }$Z_{2}^{\beta }\tilde{\partial}%
_{\beta }+Y_{2}^{b}\overset{\cdot }{\tilde{\partial}}_{b},$ \emph{we obtain
that the couple }$\left( \left[ ,\right] _{\left( \rho ,\eta \right)
TE},\left( \tilde{\rho},Id_{E}\right) \right) $\emph{\ is a Lie algebroid
structure for the vector bundle }$\left( \left( \rho ,\eta \right) TE,\left(
\rho ,\eta \right) \tau _{E},E\right) .$

\emph{Remark 2.2}\textbf{\ }In particular, if $h=Id_{M},$ then the Lie
algebroid
\begin{equation*}
\begin{array}{c}
\left( \left( \left( Id_{TM},Id_{M}\right) TE,\left( Id_{TM},Id_{M}\right)
\tau _{E},E\right) ,\left[ ,\right] _{\left( Id_{TM},Id_{M}\right)
TE},\left( \widetilde{Id_{TM}},Id_{E}\right) \right)%
\end{array}%
\end{equation*}%
is isomorphic with the usual Lie algebroid
\begin{equation*}
\begin{array}{c}
\left( \left( TE,\tau _{E},E\right) ,\left[ ,\right] _{TE},\left(
Id_{TE},Id_{E}\right) \right) .%
\end{array}%
\end{equation*}

This is a reason for which the Lie algebroid
\begin{equation*}
\begin{array}{c}
\left( \left( \left( \rho ,\eta \right) TE,\left( \rho ,\eta \right) \tau
_{E},E\right) ,\left[ ,\right] _{\left( \rho ,\eta \right) TE},\left( \tilde{%
\rho},Id_{E}\right) \right)%
\end{array}%
,
\end{equation*}%
will be called the \emph{Lie algebroid generalized tangent bundle. }(see $%
\left[ 8\right] $)

\section{$\left( \protect\rho ,\protect\eta \right) $-connections and
adapted basis}

\ \

We consider the diagram:
\begin{equation*}
\begin{array}{rcl}
E &  & \left( F,\left[ ,\right] _{F,h},\left( \rho ,\eta \right) \right) \\
\pi \downarrow &  & ~\downarrow \nu \\
M & ^{\underrightarrow{~\ \ \ \ h~\ \ \ \ }} & ~\ N%
\end{array}%
\end{equation*}%
where $\left( E,\pi ,M\right) \in \left\vert \mathbf{B}^{\mathbf{v}%
}\right\vert $ and $\left( \left( F,\nu ,N\right) ,\left[ ,\right]
_{F,h},\left( \rho ,\eta \right) \right) \in \left\vert \mathbf{GLA}%
\right\vert $.

Let
\begin{equation*}
\left( \left( \left( \rho ,\eta \right) TE,\left( \rho ,\eta \right) \tau
_{E},E\right) ,\left[ ,\right] _{\left( \rho ,\eta \right) TE},\left( \tilde{%
\rho},Id_{E}\right) \right)
\end{equation*}%
\ be the Lie algebroid generalized tangent bundle of the vector bundle $%
\left( E,\pi ,M\right) $.

We consider the $\mathbf{B}^{\mathbf{v}}$-morphism $\left( \left( \rho ,\eta
\right) \pi !,Id_{E}\right) $ given by the commutative diagram%
\begin{equation*}
\begin{array}{rcl}
\left( \rho ,\eta \right) TE & ^{\underrightarrow{~\ \left( \rho ,\eta
\right) \pi !~\ }} & \pi ^{\ast }\left( h^{\ast }F\right) \\
\left( \rho ,\eta \right) \tau _{E}\downarrow ~ &  & ~\downarrow \pi ^{\ast
}\left( h^{\ast }\nu \right) \\
E~\  & ^{\underrightarrow{~Id_{E}~}} & ~\ E%
\end{array}%
\leqno(3.1)
\end{equation*}

This is defined as:%
\begin{equation*}
\begin{array}{c}
\left( \rho ,\eta \right) \pi !\left( \left( Z^{\alpha }\tilde{\partial}%
_{\alpha }+Y^{a}\overset{\cdot }{\tilde{\partial}}_{a}\right) \left(
u_{x}\right) \right) =\left( Z^{\alpha }T_{\alpha }\right) \left(
u_{x}\right) ,%
\end{array}%
\leqno(3.2)
\end{equation*}%
for any $Z^{\alpha }\tilde{\partial}_{\alpha }+Y^{a}\overset{\cdot }{\tilde{%
\partial}}_{a}\in \Gamma \left( \left( \rho ,\eta \right) TE,\left( \rho
,\eta \right) \tau _{E},E\right) .$\medskip

Using the $\mathbf{B}^{\mathbf{v}}$-morphism $\left( \left( \rho ,\eta
\right) \pi !,Id_{E}\right) $ we obtain the \emph{tangent }$\left( \rho
,\eta \right) $\emph{-application }$\left( \left( \rho ,\eta \right) T\pi
,h\circ \pi \right) $ of $\left( \left( \rho ,\eta \right) TE,\left( \rho
,\eta \right) \tau _{E},E\right) $ source and $\left( F,\nu ,N\right) $
target.

\textbf{Definition 3.1} The kernel of the tangent $\left( \rho ,\eta \right)
$-application\ is writen
\begin{equation*}
\left( V\left( \rho ,\eta \right) TE,\left( \rho ,\eta \right) \tau
_{E},E\right)
\end{equation*}%
and it is called \emph{the vertical interior differential system}.\bigskip\
(see $\left[ 9\right] $)

We remark that the set $\left\{ \overset{\cdot }{\tilde{\partial}}_{a},~a\in
\overline{1,r}\right\} $ is a base of the $\mathcal{F}\left( E\right) $%
-module
\begin{equation*}
\left( \Gamma \left( V\left( \rho ,\eta \right) TE,\left( \rho ,\eta \right)
\tau _{E},E\right) ,+,\cdot \right) .
\end{equation*}

\textbf{Proposition 3.1} \emph{The short sequence of vector bundles}%
\begin{equation*}
\begin{array}{ccccccccc}
0 & \hookrightarrow & V\left( \rho ,\eta \right) TE & \hookrightarrow &
\left( \rho ,\eta \right) TE & ^{\underrightarrow{~\ \left( \rho ,\eta
\right) \pi !~\ }} & \pi ^{\ast }\left( h^{\ast }F\right) & ^{%
\underrightarrow{}} & 0 \\
\downarrow &  & \downarrow &  & \downarrow &  & \downarrow &  & \downarrow
\\
E & ^{\underrightarrow{~Id_{E}~}} & E & ^{\underrightarrow{~Id_{E}~}} & E &
^{\underrightarrow{~Id_{E}~}} & E & ^{\underrightarrow{~Id_{E}~}} & E%
\end{array}%
\leqno(3.3)
\end{equation*}%
\emph{is exact.}

\textbf{Definition 3.2} \textit{A }$\mathbf{Man}$-morphism $\left( \rho
,\eta \right) \Gamma $ of $\left( \rho ,\eta \right) TE$ source and $V\left(
\rho ,\eta \right) TE$ target defined by%
\begin{equation*}
\begin{array}{c}
\left( \rho ,\eta \right) \Gamma \left( Z^{\gamma }\tilde{\partial}_{\gamma
}+Y^{a}\overset{\cdot }{\tilde{\partial}}_{a}\right) \left( u_{x}\right)
=\left( Y^{a}+\left( \rho ,\eta \right) \Gamma _{\gamma }^{a}Z^{\gamma
}\right) \overset{\cdot }{\tilde{\partial}}_{a}\left( u_{x}\right) ,%
\end{array}%
\leqno(3.4)
\end{equation*}%
so that the $\mathbf{B}^{\mathbf{v}}$-morphism $\left( \left( \rho ,\eta
\right) \Gamma ,Id_{E}\right) $ is a split to the left in the previous exact
sequence, will be called $\left( \rho ,\eta \right) $\emph{-connection for
the vector bundle }$\left( E,\pi ,M\right) $.

The $\left( \rho ,Id_{M}\right) $-connection is called $\rho $\emph{%
-connection }and is denoted $\rho \Gamma $\emph{\ }and the $\left(
Id_{TM},Id_{M}\right) $-connection is called \emph{connection }and is
denoted $\Gamma $\emph{.}

\textbf{Definition 3.3 }If $\left( \rho ,\eta \right) \Gamma $ is a $\left(
\rho ,\eta \right) $-connection for the vector bundle $\left( E,\pi
,M\right) $, then the kernel of the $\mathbf{B}^{\mathbf{v}}$-morphism $%
\left( \left( \rho ,\eta \right) \Gamma ,Id_{E}\right) $\ is written $\left(
H\left( \rho ,\eta \right) TE,\left( \rho ,\eta \right) \tau _{E},E\right) $
and is called the \emph{horizontal interior differential system}. (see $%
\left[ 9\right] $)

\textbf{Definition 3.4} If $\left( E,\pi ,M\right) \in \left\vert \mathbf{B}%
^{\mathbf{v}}\right\vert $ and $\left\{ s_{a},a\in \overline{1,r}\right\} $
is a base of the $\mathcal{F}\left( M\right) $-module of sections $\left(
\Gamma \left( E,\pi ,M\right) ,+,\cdot \right) $, then the $\mathbf{B}^{%
\mathbf{v}}$-morphism $\left( \Pi ,\pi \right) $ given by the commutative
diagram%
\begin{equation*}
\begin{array}{rcl}
V\left( \rho ,\eta \right) TE & ^{\underrightarrow{~\ \Pi ~\ }} & ~\ E \\
\left( \rho ,\eta \right) \tau _{E}\downarrow ~ &  & ~\downarrow \pi \\
E~\  & ^{\underrightarrow{~~~\pi ~~}} & ~\ M%
\end{array}%
\leqno(3.5)
\end{equation*}%
is defined by
\begin{equation*}
\begin{array}{c}
\Pi \left( Y^{a}\overset{\cdot }{\tilde{\partial}}_{a}\left( u_{x}\right)
\right) =Y^{a}\left( u_{x}\right) s_{a}\left( x\right) .%
\end{array}%
\leqno(3.6)
\end{equation*}

\textbf{Theorem 3.1 }(see $\left[ 8\right] $) \emph{If }$\left( \rho ,\eta
\right) \Gamma $\emph{\ is a }$\left( \rho ,\eta \right) $\emph{-connection
for the vector bundle }$\left( E,\pi ,M\right) ,$\emph{\ then its components
satisfy the law of transformation }%
\begin{equation*}
\begin{array}{c}
\left( \rho ,\eta \right) \Gamma _{\gamma
%TCIMACRO{\U{b4}}%
%BeginExpansion
{\acute{}}%
%EndExpansion
}^{a%
%TCIMACRO{\U{b4}}%
%BeginExpansion
{\acute{}}%
%EndExpansion
}{=}M_{a}^{a%
%TCIMACRO{\U{b4}}%
%BeginExpansion
{\acute{}}%
%EndExpansion
}{\circ }\pi \!\!\left[ \rho _{\gamma }^{k}{\circ }h{\circ }\pi \frac{%
\partial M_{b%
%TCIMACRO{\U{b4}}%
%BeginExpansion
{\acute{}}%
%EndExpansion
}^{a}\circ \pi }{\partial x^{k}}y^{b%
%TCIMACRO{\U{b4}}%
%BeginExpansion
{\acute{}}%
%EndExpansion
}{+}\left( \rho ,\eta \right) \!\Gamma _{\gamma }^{a}\right] \!\!\Lambda
_{\gamma
%TCIMACRO{\U{b4}}%
%BeginExpansion
{\acute{}}%
%EndExpansion
}^{\gamma }{\circ }h{\circ }\pi .%
\end{array}%
\leqno(3.7)
\end{equation*}

\emph{In the particular case of Lie algebroids, }$\left( \eta ,h\right)
=\left( Id_{M},Id_{M}\right) ,$\emph{\ the relations }$\left( 3.7\right) $%
\emph{\ become}%
\begin{equation*}
\begin{array}{c}
\rho \Gamma _{\gamma
%TCIMACRO{\U{b4}}%
%BeginExpansion
{\acute{}}%
%EndExpansion
}^{a%
%TCIMACRO{\U{b4}}%
%BeginExpansion
{\acute{}}%
%EndExpansion
}=M_{a}^{a%
%TCIMACRO{\U{b4}}%
%BeginExpansion
{\acute{}}%
%EndExpansion
}\circ \pi \left[ \rho _{\gamma }^{k}\circ \pi \frac{\partial M_{b%
%TCIMACRO{\U{b4}}%
%BeginExpansion
{\acute{}}%
%EndExpansion
}^{a}\circ \pi }{\partial x^{k}}y^{b%
%TCIMACRO{\U{b4}}%
%BeginExpansion
{\acute{}}%
%EndExpansion
}+\rho \Gamma _{\gamma }^{a}\right] \Lambda _{\gamma
%TCIMACRO{\U{b4}}%
%BeginExpansion
{\acute{}}%
%EndExpansion
}^{\gamma }\circ \pi .%
\end{array}%
\leqno(3.7^{\prime })
\end{equation*}

\emph{In the classical case, }$\left( \rho ,\eta ,h\right) =\left(
Id_{TM},Id_{M},Id_{M}\right) ,$\emph{\ the relations }$\left( 3.7^{\prime
}\right) $\emph{\ become}%
\begin{equation*}
\begin{array}{c}
\Gamma _{k%
%TCIMACRO{\U{b4}}%
%BeginExpansion
{\acute{}}%
%EndExpansion
}^{i%
%TCIMACRO{\U{b4}}%
%BeginExpansion
{\acute{}}%
%EndExpansion
}=\frac{\partial x^{i%
%TCIMACRO{\U{b4}}%
%BeginExpansion
{\acute{}}%
%EndExpansion
}}{\partial x^{i}}\circ \tau _{M}\left[ \frac{\partial }{\partial x^{k}}%
\left( \frac{\partial x^{i}}{\partial x^{j%
%TCIMACRO{\U{b4}}%
%BeginExpansion
{\acute{}}%
%EndExpansion
}}\circ \tau _{M}\right) y^{j%
%TCIMACRO{\U{b4}}%
%BeginExpansion
{\acute{}}%
%EndExpansion
}+\Gamma _{k}^{i}\right] \frac{\partial x^{k}}{\partial x^{k%
%TCIMACRO{\U{b4}}%
%BeginExpansion
{\acute{}}%
%EndExpansion
}}\circ \tau _{M}.%
\end{array}%
\leqno(3.7^{\prime \prime })
\end{equation*}

\emph{Remark 3.1} \emph{If we have a set of real local functions }$\left(
\rho ,\eta \right) \Gamma _{\gamma }^{a}$\emph{\ which satisfies the
relations of passing }$\left( 3.7\right) ,$\emph{\ then we have a }$\left(
\rho ,\eta \right) $\emph{-connection }$\left( \rho ,\eta \right) \Gamma $%
\emph{\ for the vector bundle} $\left( E,\pi ,M\right) $

\textbf{Example 3.1 }If $\Gamma $ is an Ehresmann connection for the vector
bundle $\left( E,\pi ,M\right) $ on components $\Gamma _{k}^{a},$ then the
differentiable real local functions $\left( \rho ,\eta \right) \Gamma
_{\gamma }^{a}=\left( \rho _{\gamma }^{k}\circ h\circ \pi \right) \Gamma
_{k}^{a}$ are the components of a $\left( \rho ,\eta \right) $-connection $%
\left( \rho ,\eta \right) \Gamma $ for the vector bundle $\left( E,\pi
,M\right) .$ This $\left( \rho ,\eta \right) $-connection will be called the
$\left( \rho ,\eta \right) $\emph{-connection associated to the connection }$%
\Gamma .$

We put the problem of finding a base for the $\mathcal{F}\left( E\right) $%
-module
\begin{equation*}
\left( \Gamma \left( H\left( \rho ,\eta \right) TE,\left( \rho ,\eta \right)
\tau _{E},E\right) ,+,\cdot \right)
\end{equation*}%
of the type\textbf{\ }
\begin{equation*}
\begin{array}[t]{l}
\tilde{\delta}_{\alpha }=Z_{\alpha }^{\beta }\tilde{\partial}_{\beta
}+Y_{\alpha }^{a}\overset{\cdot }{\tilde{\partial}}_{a},\alpha \in \overline{%
1,r}%
\end{array}%
\end{equation*}%
which satisfies the following conditions:
\begin{equation*}
\begin{array}{rcl}
\displaystyle\Gamma \left( \left( \rho ,\eta \right) \pi !,Id_{E}\right)
\left( \tilde{\delta}_{\alpha }\right) & = & T_{\alpha }\vspace*{2mm}, \\
\displaystyle\Gamma \left( \left( \rho ,\eta \right) \Gamma ,Id_{E}\right)
\left( \tilde{\delta}_{\alpha }\right) & = & 0.%
\end{array}%
\leqno(3.8)
\end{equation*}

Then we obtain the sections
\begin{equation*}
\begin{array}[t]{l}
\frac{\delta }{\delta \tilde{z}^{\alpha }}=\tilde{\partial}_{\alpha }-\left(
\rho ,\eta \right) \Gamma _{\alpha }^{a}\overset{\cdot }{\tilde{\partial}}%
_{a}=T_{\alpha }\oplus \left( \left( \rho _{\alpha }^{i}\circ h\circ \pi
\right) \partial _{i}-\left( \rho ,\eta \right) \Gamma _{\alpha }^{a}\dot{%
\partial}_{a}\right) .%
\end{array}%
\leqno(3.9)
\end{equation*}%
such that their law of change is a tensorial law under a change of vector
fiber charts.

The base $\left( \tilde{\delta}_{\alpha },\overset{\cdot }{\tilde{\partial}}%
_{a}\right) $ will be called the \emph{adapted }$\left( \rho ,\eta \right) $%
\emph{-base.}

\emph{Remark 3.2 The following equality holds good}%
\begin{equation*}
\begin{array}{l}
\Gamma \left( \tilde{\rho},Id_{E}\right) \left( \tilde{\delta}_{\alpha
}\right) =\left( \rho _{\alpha }^{i}\circ h\circ \pi \right) \partial
_{i}-\left( \rho ,\eta \right) \Gamma _{\alpha }^{a}\dot{\partial}_{a}.%
\end{array}%
\leqno(3.10)
\end{equation*}

Moreover, if $\left( \rho ,\eta \right) \Gamma $ is the $\left( \rho ,\eta
\right) $-connection associated to a connection $\Gamma $, then we obtain
\begin{equation*}
\begin{array}{l}
\Gamma \left( \tilde{\rho},Id_{E}\right) \left( \tilde{\delta}_{\alpha
}\right) =\left( \rho _{\alpha }^{i}\circ h\circ \pi \right) \delta _{i},%
\end{array}%
\leqno(3.11)
\end{equation*}%
where $\left( \delta _{i},\dot{\partial}_{a}\right) $ is the adapted base
for the $\mathcal{F}\left( E\right) $-module $\left( \Gamma \left( TE,\tau
_{E},E\right) ,+,\cdot \right) .$

Let $\left( d\tilde{z}^{\alpha },d\tilde{y}^{b}\right) $ be the natural dual
$\left( \rho ,\eta \right) $-base of natural $\left( \rho ,\eta \right) $%
-base $\left( \displaystyle\tilde{\partial}_{\alpha },\displaystyle\overset{%
\cdot }{\tilde{\partial}}_{a}\right) .$

This is determined by the equations
\begin{equation*}
\begin{array}{c}
\left\{
\begin{array}{cc}
\displaystyle\left\langle d\tilde{z}^{\alpha },\tilde{\partial}_{\beta
}\right\rangle =\delta _{\beta }^{\alpha }, & \displaystyle\left\langle d%
\tilde{z}^{\alpha },\overset{\cdot }{\tilde{\partial}}_{a}\right\rangle =0,%
\vspace*{2mm} \\
\displaystyle\left\langle d\tilde{y}^{a},\tilde{\partial}_{\beta
}\right\rangle =0, & \displaystyle\left\langle d\tilde{y}^{a},\overset{\cdot
}{\tilde{\partial}}_{b}\right\rangle =\delta _{b}^{a}.%
\end{array}%
\right.%
\end{array}%
\end{equation*}

We consider the problem of finding a base for the $\mathcal{F}\left(
E\right) $-module
\begin{equation*}
\left( \Gamma \left( \left( V\left( \rho ,\eta \right) TE\right) ^{\ast
},\left( \left( \rho ,\eta \right) \tau _{E}\right) ^{\ast },E\right)
,+,\cdot \right)
\end{equation*}%
of the type
\begin{equation*}
\begin{array}{c}
\delta \tilde{y}^{a}=\theta _{\alpha }^{a}d\tilde{z}^{\alpha }+\omega
_{b}^{a}d\tilde{y}^{b},~a\in \overline{1,n}%
\end{array}%
\end{equation*}%
which satisfies the following conditions:
\begin{equation*}
\begin{array}{c}
\left\langle \delta \tilde{y}^{a},\overset{\cdot }{\tilde{\partial}}%
_{a}\right\rangle =1\wedge \left\langle \delta \tilde{y}^{a},\tilde{\delta}%
_{\alpha }\right\rangle =0.%
\end{array}%
\leqno(3.12)
\end{equation*}

We obtain the sections
\begin{equation*}
\begin{array}{l}
\delta \tilde{y}^{a}=\left( \rho ,\eta \right) \Gamma _{\alpha }^{a}d\tilde{z%
}^{\alpha }+d\tilde{y}^{a},a\in \overline{1,n}.%
\end{array}%
\leqno(3.13)
\end{equation*}%
such that their changing rule is tensorial under a change of vector fiber
charts. The base $\left( d\tilde{z}^{\alpha },\delta \tilde{y}^{a}\right) $
will be called the \emph{adapted dual }$\left( \rho ,\eta \right) $\emph{%
-base.}

\section{Distinguished linear $\left( \protect\rho ,\protect\eta \right) $%
-connections}

\ \ \

We consider the following diagram:
\begin{equation*}
\begin{array}{rcl}
E &  & \left( F,\left[ ,\right] _{F,h},\left( \rho ,\eta \right) \right) \\
\pi \downarrow &  & ~\downarrow \nu \\
M & ^{\underrightarrow{~\ \ \ \ h~\ \ \ \ }} & ~\ N%
\end{array}%
\leqno(4.1)
\end{equation*}%
where $\left( E,\pi ,M\right) \in \left\vert \mathbf{B}^{\mathbf{v}%
}\right\vert $ and $\left( \left( F,\nu ,N\right) ,\left[ ,\right]
_{F,h},\left( \rho ,\eta \right) \right) $ is a generalized Lie algebroid.
Let $\left( \rho ,\eta \right) \Gamma $ be a $\left( \rho ,\eta \right) $%
-connection for the vector bundle $\left( E,\pi ,M\right) .$

Let
\begin{equation*}
\left( \mathcal{T}~_{q,s}^{p,r}\left( \left( \rho ,\eta \right) TE,\left(
\rho ,\eta \right) \tau _{E},E\right) ,+,\cdot \right)
\end{equation*}%
be the $\mathcal{F}\left( E\right) $-module of tensor fields by $\left(
_{q,s}^{p,r}\right) $-type from the generalized tangent bundle
\begin{equation*}
\left( H\left( \rho ,\eta \right) TE,\left( \rho ,\eta \right) \tau
_{E},E\right) \oplus \left( V\left( \rho ,\eta \right) TE,\left( \rho ,\eta
\right) \tau _{E},E\right) .
\end{equation*}

An arbitrarily tensor field $T$ is written as
\begin{equation*}
\begin{array}{c}
T=T_{\beta _{1}...\beta _{q}b_{1}...b_{s}}^{\alpha _{1}...\alpha
_{p}a_{1}...a_{r}}\tilde{\delta}_{\alpha _{1}}\otimes ...\otimes \tilde{%
\delta}_{\alpha _{p}}\otimes d\tilde{z}^{\beta _{1}}\otimes ...\otimes d%
\tilde{z}^{\beta _{q}}\otimes \\
\overset{\cdot }{\tilde{\partial}}_{a_{1}}\otimes ...\otimes \overset{\cdot }%
{\tilde{\partial}}_{a_{r}}\otimes \delta \tilde{y}^{b_{1}}\otimes ...\otimes
\delta \tilde{y}^{b_{s}}.%
\end{array}%
\end{equation*}

Let
\begin{equation*}
\left( \mathcal{T}~\left( \left( \rho ,\eta \right) TE,\left( \rho ,\eta
\right) \tau _{E},E\right) ,+,\cdot ,\otimes \right)
\end{equation*}%
be the tensor fields algebra of generalized tangent bundle $\left( \left(
\rho ,\eta \right) TE,\left( \rho ,\eta \right) \tau _{E},E\right) $.

\textbf{Definition 4.1 }Let
\begin{equation*}
\begin{array}{l}
\left( X,T\right) ^{\ \underrightarrow{\left( \rho ,\eta \right) D}\,}%
\vspace*{1mm}\left( \rho ,\eta \right) D_{X}T%
\end{array}%
\end{equation*}%
be a covariant $\left( \rho ,\eta \right) $-derivative for the tensor
algebra
\begin{equation*}
\left( \mathcal{T}~\left( \left( \rho ,\eta \right) TE,\left( \rho ,\eta
\right) \tau _{E},E\right) ,+,\cdot ,\otimes \right)
\end{equation*}%
of the generalized tangent bundle
\begin{equation*}
\left( \left( \rho ,\eta \right) TE,\left( \rho ,\eta \right) \tau
_{E},E\right)
\end{equation*}%
which preserves the horizontal and vertical interior differential systems by
parallelism. (see $\left[ 9\right] )$

The real local functions
\begin{equation*}
\left( \left( \rho ,\eta \right) H_{\beta \gamma }^{\alpha },\left( \rho
,\eta \right) H_{b\gamma }^{a},\left( \rho ,\eta \right) V_{\beta c}^{\alpha
},\left( \rho ,\eta \right) V_{bc}^{a}\right)
\end{equation*}%
defined by the following equalities:
\begin{equation*}
\begin{array}{ll}
\left( \rho ,\eta \right) D_{\tilde{\delta}_{\gamma }}\tilde{\delta}_{\beta
}=\left( \rho ,\eta \right) H_{\beta \gamma }^{\alpha }\tilde{\delta}%
_{\alpha }, & \left( \rho ,\eta \right) D_{\tilde{\delta}_{\gamma }}\overset{%
\cdot }{\tilde{\partial}}_{b}=\left( \rho ,\eta \right) H_{b\gamma }^{a}%
\overset{\cdot }{\tilde{\partial}}_{a} \\
\left( \rho ,\eta \right) D_{\overset{\cdot }{\tilde{\partial}}_{c}}\tilde{%
\delta}_{\beta }=\left( \rho ,\eta \right) V_{\beta c}^{\alpha }\tilde{\delta%
}_{\alpha }, & \left( \rho ,\eta \right) D_{\overset{\cdot }{\tilde{\partial}%
}_{c}}\overset{\cdot }{\tilde{\partial}}_{b}=\left( \rho ,\eta \right)
V_{bc}^{a}\overset{\cdot }{\tilde{\partial}}_{a}%
\end{array}%
\leqno(4.2)
\end{equation*}%
are the components of a linear $\left( \rho ,\eta \right) $-connection $%
\left( \left( \rho ,\eta \right) H,\left( \rho ,\eta \right) V\right) $ for
the generalized tangent bundle $\left( \left( \rho ,\eta \right) TE,\left(
\rho ,\eta \right) \tau _{E},E\right) $ which will be called the \emph{%
distinguished linear }$\left( \rho ,\eta \right) $\emph{-connection.}

In the particular case of Lie algebroids, $h=Id_{M}=\eta ,$ we obtain the
\emph{distinguished linear }$\rho $\emph{-connection.} The components of a
distinguished linear $\rho $-connection $\left( \rho H,\rho V\right) $ will
be denoted
\begin{equation*}
\left( \rho H_{\beta \gamma }^{\alpha },\rho H_{b\gamma }^{a},\rho V_{\beta
c}^{\alpha },\rho V_{bc}^{a}\right) .
\end{equation*}

In addition, if $\rho =Id_{TM},$ then we obtain the classical \emph{%
distinguished linear connection. }The components of a distinguished linear
connection $\left( H,V\right) $ will be denoted
\begin{equation*}
\left( H_{jk}^{i},H_{bk}^{a},V_{jc}^{i},V_{bc}^{a}\right) .
\end{equation*}

\textbf{Theorem 4.1 }\emph{If }$((\rho ,\eta )H,(\rho ,\eta )V)$ \emph{is a
distinguished linear} $(\rho ,\eta )$-\emph{connection for the generalized
tangent bundle }$\left( \left( \rho ,\eta \right) TE,\left( \rho ,\eta
\right) \tau _{E},E\right) $\emph{, then its components satisfy the change
relations: }

\begin{equation*}
\begin{array}{ll}
\left( \rho ,\eta \right) H_{\beta
%TCIMACRO{\U{b4}}%
%BeginExpansion
{\acute{}}%
%EndExpansion
\gamma
%TCIMACRO{\U{b4}}%
%BeginExpansion
{\acute{}}%
%EndExpansion
}^{\alpha
%TCIMACRO{\U{b4}}%
%BeginExpansion
{\acute{}}%
%EndExpansion
}\!\! & =\Lambda _{\alpha }^{\alpha
%TCIMACRO{\U{b4}}%
%BeginExpansion
{\acute{}}%
%EndExpansion
}\circ h\circ \pi \cdot \left[ \Gamma \left( \tilde{\rho},Id_{E}\right)
\left( \tilde{\delta}_{\gamma }\right) \left( \Lambda _{\beta
%TCIMACRO{\U{b4}}%
%BeginExpansion
{\acute{}}%
%EndExpansion
}^{\alpha }\circ h\circ \pi \right) +\right. \vspace*{1mm} \\
& +\left. \left( \rho ,\eta \right) H_{\beta \gamma }^{\alpha }\cdot \Lambda
_{\beta
%TCIMACRO{\U{b4}}%
%BeginExpansion
{\acute{}}%
%EndExpansion
}^{\beta }\circ h\circ \pi \right] \cdot \Lambda _{\gamma
%TCIMACRO{\U{b4}}%
%BeginExpansion
{\acute{}}%
%EndExpansion
}^{\gamma }\circ h\circ \pi ,\vspace*{2mm} \\
\left( \rho ,\eta \right) H_{b%
%TCIMACRO{\U{b4}}%
%BeginExpansion
{\acute{}}%
%EndExpansion
\gamma
%TCIMACRO{\U{b4}}%
%BeginExpansion
{\acute{}}%
%EndExpansion
}^{a%
%TCIMACRO{\U{b4}}%
%BeginExpansion
{\acute{}}%
%EndExpansion
}\!\! & =M_{a}^{a%
%TCIMACRO{\U{b4}}%
%BeginExpansion
{\acute{}}%
%EndExpansion
}\circ \pi \cdot \left[ \Gamma \left( \tilde{\rho},Id_{E}\right) \left(
\tilde{\delta}_{\gamma }\right) \left( M_{b%
%TCIMACRO{\U{b4}}%
%BeginExpansion
{\acute{}}%
%EndExpansion
}^{a}\circ \pi \right) +\right. \vspace*{1mm} \\
& \left. +\left( \rho ,\eta \right) H_{b\gamma }^{a}\cdot M_{b%
%TCIMACRO{\U{b4}}%
%BeginExpansion
{\acute{}}%
%EndExpansion
}^{b}\circ \pi \right] \cdot \Lambda _{\gamma
%TCIMACRO{\U{b4}}%
%BeginExpansion
{\acute{}}%
%EndExpansion
}^{\gamma }\circ h\circ \pi ,\vspace*{2mm} \\
\left( \rho ,\eta \right) V_{\beta
%TCIMACRO{\U{b4}}%
%BeginExpansion
{\acute{}}%
%EndExpansion
c%
%TCIMACRO{\U{b4}}%
%BeginExpansion
{\acute{}}%
%EndExpansion
}^{\alpha
%TCIMACRO{\U{b4}}%
%BeginExpansion
{\acute{}}%
%EndExpansion
}\!\! & =\Lambda _{\alpha
%TCIMACRO{\U{b4}}%
%BeginExpansion
{\acute{}}%
%EndExpansion
}^{\alpha }\circ h\circ \pi \cdot \left( \rho ,\eta \right) V_{\beta
c}^{\alpha }\cdot \Lambda _{\beta
%TCIMACRO{\U{b4}}%
%BeginExpansion
{\acute{}}%
%EndExpansion
}^{\beta }\circ h\circ \pi \cdot M_{c%
%TCIMACRO{\U{b4}}%
%BeginExpansion
{\acute{}}%
%EndExpansion
}^{c}\circ \pi ,\vspace*{2mm} \\
\left( \rho ,\eta \right) V_{b%
%TCIMACRO{\U{b4}}%
%BeginExpansion
{\acute{}}%
%EndExpansion
c%
%TCIMACRO{\U{b4}}%
%BeginExpansion
{\acute{}}%
%EndExpansion
}^{a%
%TCIMACRO{\U{b4}}%
%BeginExpansion
{\acute{}}%
%EndExpansion
}\!\! & =M_{a}^{a%
%TCIMACRO{\U{b4}}%
%BeginExpansion
{\acute{}}%
%EndExpansion
}\circ \pi \cdot \left( \rho ,\eta \right) V_{bc}^{a}\cdot M_{b%
%TCIMACRO{\U{b4}}%
%BeginExpansion
{\acute{}}%
%EndExpansion
}^{b}\circ \pi \cdot M_{c%
%TCIMACRO{\U{b4}}%
%BeginExpansion
{\acute{}}%
%EndExpansion
}^{c}\circ \pi .%
\end{array}%
\leqno(4.3)
\end{equation*}

\textbf{Corollary 4.1 }\emph{In the particular case of Lie algebroids, }$%
\left( \eta ,h\right) =\left( Id_{M},Id_{M}\right) ,$\emph{\ we obtain }%
\begin{equation*}
\begin{array}{ll}
\rho H_{\beta
%TCIMACRO{\U{b4}}%
%BeginExpansion
{\acute{}}%
%EndExpansion
\gamma
%TCIMACRO{\U{b4}}%
%BeginExpansion
{\acute{}}%
%EndExpansion
}^{\alpha
%TCIMACRO{\U{b4}}%
%BeginExpansion
{\acute{}}%
%EndExpansion
}\!\! & =\Lambda _{\alpha }^{\alpha
%TCIMACRO{\U{b4}}%
%BeginExpansion
{\acute{}}%
%EndExpansion
}\circ \pi \cdot \left[ \Gamma \left( \tilde{\rho},Id_{E}\right) \left(
\tilde{\delta}_{\gamma }\right) \left( \Lambda _{\beta
%TCIMACRO{\U{b4}}%
%BeginExpansion
{\acute{}}%
%EndExpansion
}^{\alpha }\circ \pi \right) +\rho H_{\beta \gamma }^{\alpha }\cdot \Lambda
_{\beta
%TCIMACRO{\U{b4}}%
%BeginExpansion
{\acute{}}%
%EndExpansion
}^{\beta }\circ \pi \right] \cdot \Lambda _{\gamma
%TCIMACRO{\U{b4}}%
%BeginExpansion
{\acute{}}%
%EndExpansion
}^{\gamma }\circ \pi \\
\rho H_{b%
%TCIMACRO{\U{b4}}%
%BeginExpansion
{\acute{}}%
%EndExpansion
\gamma
%TCIMACRO{\U{b4}}%
%BeginExpansion
{\acute{}}%
%EndExpansion
}^{a%
%TCIMACRO{\U{b4}}%
%BeginExpansion
{\acute{}}%
%EndExpansion
}\!\! & =M_{a}^{a%
%TCIMACRO{\U{b4}}%
%BeginExpansion
{\acute{}}%
%EndExpansion
}\circ \pi \cdot \left[ \Gamma \left( \tilde{\rho},Id_{E}\right) \left(
\tilde{\delta}_{\gamma }\right) \left( M_{b%
%TCIMACRO{\U{b4}}%
%BeginExpansion
{\acute{}}%
%EndExpansion
}^{a}\circ \pi \right) +\rho H_{b\gamma }^{a}\cdot M_{b%
%TCIMACRO{\U{b4}}%
%BeginExpansion
{\acute{}}%
%EndExpansion
}^{b}\circ \pi \right] \cdot \Lambda _{\gamma
%TCIMACRO{\U{b4}}%
%BeginExpansion
{\acute{}}%
%EndExpansion
}^{\gamma }\circ \pi , \\
\rho V_{\beta
%TCIMACRO{\U{b4}}%
%BeginExpansion
{\acute{}}%
%EndExpansion
c%
%TCIMACRO{\U{b4}}%
%BeginExpansion
{\acute{}}%
%EndExpansion
}^{\alpha
%TCIMACRO{\U{b4}}%
%BeginExpansion
{\acute{}}%
%EndExpansion
}\!\! & =\Lambda _{\alpha
%TCIMACRO{\U{b4}}%
%BeginExpansion
{\acute{}}%
%EndExpansion
}^{\alpha }\circ \pi \cdot \rho V_{\beta c}^{\alpha }\cdot \Lambda _{\beta
%TCIMACRO{\U{b4}}%
%BeginExpansion
{\acute{}}%
%EndExpansion
}^{\beta }\circ \pi \cdot M_{c%
%TCIMACRO{\U{b4}}%
%BeginExpansion
{\acute{}}%
%EndExpansion
}^{c}\circ \pi ,\vspace*{2mm} \\
\rho V_{b%
%TCIMACRO{\U{b4}}%
%BeginExpansion
{\acute{}}%
%EndExpansion
c%
%TCIMACRO{\U{b4}}%
%BeginExpansion
{\acute{}}%
%EndExpansion
}^{a%
%TCIMACRO{\U{b4}}%
%BeginExpansion
{\acute{}}%
%EndExpansion
}\!\! & =M_{a}^{a%
%TCIMACRO{\U{b4}}%
%BeginExpansion
{\acute{}}%
%EndExpansion
}\circ \pi \cdot \rho V_{bc}^{a}\cdot M_{b%
%TCIMACRO{\U{b4}}%
%BeginExpansion
{\acute{}}%
%EndExpansion
}^{b}\circ \pi \cdot M_{c%
%TCIMACRO{\U{b4}}%
%BeginExpansion
{\acute{}}%
%EndExpansion
}^{c}\circ \pi .%
\end{array}%
\leqno(4.3^{\prime })
\end{equation*}

\emph{In the classical case, }$\left( \rho ,\eta ,h\right) =\left(
Id_{TM},Id_{M},Id_{M}\right) ,$\emph{\ we obtain that the components of a
distinguished linear connection }$\left( H,V\right) $\emph{\ verify the
change relations:}
\begin{equation*}
\begin{array}{cl}
H_{j%
%TCIMACRO{\U{b4}}%
%BeginExpansion
{\acute{}}%
%EndExpansion
k%
%TCIMACRO{\U{b4}}%
%BeginExpansion
{\acute{}}%
%EndExpansion
}^{i%
%TCIMACRO{\U{b4}}%
%BeginExpansion
{\acute{}}%
%EndExpansion
} & =\frac{\partial x^{i%
%TCIMACRO{\U{b4}}%
%BeginExpansion
{\acute{}}%
%EndExpansion
}}{\partial x^{i}}\circ \pi \cdot \left[ \frac{\delta }{\delta x^{k}}\left(
\frac{\partial x^{i}}{\partial x^{j%
%TCIMACRO{\U{b4}}%
%BeginExpansion
{\acute{}}%
%EndExpansion
}}\circ \pi \right) +H_{jk}^{i}\cdot \frac{\partial x^{j}}{\partial x^{j%
%TCIMACRO{\U{b4}}%
%BeginExpansion
{\acute{}}%
%EndExpansion
}}\circ \pi \right] \cdot \frac{\partial x^{k}}{\partial x^{k%
%TCIMACRO{\U{b4}}%
%BeginExpansion
{\acute{}}%
%EndExpansion
}}\circ \pi ,\vspace*{2mm} \\
H_{b%
%TCIMACRO{\U{b4}}%
%BeginExpansion
{\acute{}}%
%EndExpansion
k%
%TCIMACRO{\U{b4}}%
%BeginExpansion
{\acute{}}%
%EndExpansion
}^{a%
%TCIMACRO{\U{b4}}%
%BeginExpansion
{\acute{}}%
%EndExpansion
} & =M_{a}^{a%
%TCIMACRO{\U{b4}}%
%BeginExpansion
{\acute{}}%
%EndExpansion
}\circ \pi \cdot \left[ \frac{\delta }{\delta x^{k}}\left( M_{b%
%TCIMACRO{\U{b4}}%
%BeginExpansion
{\acute{}}%
%EndExpansion
}^{a}\circ \pi \right) +H_{bk}^{a}\cdot M_{b%
%TCIMACRO{\U{b4}}%
%BeginExpansion
{\acute{}}%
%EndExpansion
}^{b}\circ \pi \right] \cdot \frac{\partial x^{k}}{\partial x^{k%
%TCIMACRO{\U{b4}}%
%BeginExpansion
{\acute{}}%
%EndExpansion
}}\circ \pi ,\vspace*{2mm} \\
V_{j%
%TCIMACRO{\U{b4}}%
%BeginExpansion
{\acute{}}%
%EndExpansion
c%
%TCIMACRO{\U{b4}}%
%BeginExpansion
{\acute{}}%
%EndExpansion
}^{i%
%TCIMACRO{\U{b4}}%
%BeginExpansion
{\acute{}}%
%EndExpansion
} & =\frac{\partial x^{i%
%TCIMACRO{\U{b4}}%
%BeginExpansion
{\acute{}}%
%EndExpansion
}}{\partial x^{i}}\circ \pi \cdot V_{jc}^{i}\frac{\partial x^{j}}{\partial
x^{j%
%TCIMACRO{\U{b4}}%
%BeginExpansion
{\acute{}}%
%EndExpansion
}}\circ \pi \cdot M_{c%
%TCIMACRO{\U{b4}}%
%BeginExpansion
{\acute{}}%
%EndExpansion
}^{c}\circ \pi ,\vspace*{3mm} \\
V_{b%
%TCIMACRO{\U{b4}}%
%BeginExpansion
{\acute{}}%
%EndExpansion
c%
%TCIMACRO{\U{b4}}%
%BeginExpansion
{\acute{}}%
%EndExpansion
}^{a%
%TCIMACRO{\U{b4}}%
%BeginExpansion
{\acute{}}%
%EndExpansion
} & =M_{a}^{a%
%TCIMACRO{\U{b4}}%
%BeginExpansion
{\acute{}}%
%EndExpansion
}\circ \pi \cdot V_{bc}^{a}\cdot M_{b%
%TCIMACRO{\U{b4}}%
%BeginExpansion
{\acute{}}%
%EndExpansion
}^{b}\circ \pi \cdot M_{c%
%TCIMACRO{\U{b4}}%
%BeginExpansion
{\acute{}}%
%EndExpansion
}^{c}\circ \pi .%
\end{array}%
\leqno(4.3^{\prime \prime })
\end{equation*}

\textbf{Example 4.1} The local real functions%
\begin{equation*}
\begin{array}[b]{c}
\left( \frac{\partial \left( \rho ,\eta \right) \Gamma _{\gamma }^{a}}{%
\partial y^{b}},\frac{\partial \left( \rho ,\eta \right) \Gamma _{\gamma
}^{a}}{\partial y^{b}},0,0\right)%
\end{array}%
\leqno(4.4)
\end{equation*}%
are the components of a distinguished linear $\left( \rho ,\eta \right) $%
\textit{-}connection for the generalized tangent bundle $\left( \left( \rho
,\eta \right) TE,\left( \rho ,\eta \right) \tau _{E},E\right) ,$ which will
by called the \emph{Berwald linear }$\left( \rho ,\eta \right) $\emph{%
-connection.}

The Berwald linear $(Id_{TM},Id_{M})$-connection are the usual \emph{%
Ber\-wald linear connection.}

\textbf{Theorem 4.2} \emph{If the generalized tangent bundle} $\!(\!(\rho
,\!\eta )T\!E,\!(\rho ,\!\eta )\tau _{E},\!E\!)$ \emph{is endowed with a
distinguished linear} $\!(\rho ,\!\eta )$\emph{-connection} $((\rho ,\eta
)H,(\rho ,\eta )V),$ \emph{then for any}
\begin{equation*}
\begin{array}[b]{c}
X=Z^{\alpha }\tilde{\delta}_{\alpha }+Y^{a}\overset{\cdot }{\tilde{\partial}}%
_{a}\in \Gamma (\!(\rho ,\eta )TE,\!(\rho ,\!\eta )\tau _{E},\!E)%
\end{array}%
\end{equation*}%
\emph{and for any}
\begin{equation*}
T\in \mathcal{T}_{qs}^{pr}\!(\!(\rho ,\eta )TE,\!(\rho ,\eta )\tau _{E},\!E),
\end{equation*}%
\emph{we obtain the formula:}
\begin{equation*}
\begin{array}{l}
\left( \rho ,\eta \right) D_{X}\left( T_{\beta _{1}...\beta
_{q}b_{1}...b_{s}}^{\alpha _{1}...\alpha _{p}a_{1}...a_{r}}\tilde{\delta}%
_{\alpha _{1}}\otimes ...\otimes \tilde{\delta}_{\alpha _{p}}\otimes d\tilde{%
z}^{\beta _{1}}\otimes ...\otimes \right. \vspace*{1mm} \\
\hspace*{9mm}\left. \otimes d\tilde{z}^{\beta _{q}}\otimes \overset{\cdot }{%
\tilde{\partial}}_{a_{1}}\otimes ...\otimes \overset{\cdot }{\tilde{\partial}%
}_{a_{r}}\otimes \delta \tilde{y}^{b_{1}}\otimes ...\otimes \delta \tilde{y}%
^{b_{s}}\right) =\vspace*{1mm} \\
\hspace*{9mm}=Z^{\gamma }T_{\beta _{1}...\beta _{q}b_{1}...b_{s}\mid \gamma
}^{\alpha _{1}...\alpha _{p}a_{1}...a_{r}}\tilde{\delta}_{\alpha
_{1}}\otimes ...\otimes \tilde{\delta}_{\alpha _{p}}\otimes d\tilde{z}%
^{\beta _{1}}\otimes ...\otimes d\tilde{z}^{\beta _{q}}\otimes \overset{%
\cdot }{\tilde{\partial}}_{a_{1}}\otimes ...\otimes \vspace*{1mm} \\
\hspace*{9mm}\otimes \overset{\cdot }{\tilde{\partial}}_{a_{r}}\otimes
\delta \tilde{y}^{b_{1}}\otimes ...\otimes \delta \tilde{y}%
^{b_{s}}+Y^{c}T_{\beta _{1}...\beta _{q}b_{1}...b_{s}}^{\alpha _{1}...\alpha
_{p}a_{1}...a_{r}}\mid _{c}\tilde{\delta}_{\alpha _{1}}\otimes ...\otimes
\vspace*{1mm} \\
\hspace*{9mm}\otimes \tilde{\delta}_{\alpha _{p}}\otimes d\tilde{z}^{\beta
_{1}}\otimes ...\otimes d\tilde{z}^{\beta _{q}}\otimes \overset{\cdot }{%
\tilde{\partial}}_{a_{1}}\otimes ...\otimes \overset{\cdot }{\tilde{\partial}%
}_{a_{r}}\otimes \delta \tilde{y}^{b_{1}}\otimes ...\otimes \delta \tilde{y}%
^{b_{s}},%
\end{array}%
\leqno(4.5)
\end{equation*}%
\emph{where }%
\begin{equation*}
\begin{array}{l}
T_{\beta _{1}...\beta _{q}b_{1}...b_{s}\mid \gamma }^{\alpha _{1}...\alpha
_{p}a_{1}...a_{r}}=\vspace*{2mm}\Gamma \left( \tilde{\rho},Id_{E}\right)
\left( \tilde{\delta}_{\gamma }\right) T_{\beta _{1}...\beta
_{q}b_{1}...b_{s}}^{\alpha _{1}...\alpha _{p}a_{1}...a_{r}} \\
\hspace*{8mm}+\left( \rho ,\eta \right) H_{\alpha \gamma }^{\alpha
_{1}}T_{\beta _{1}...\beta _{q}b_{1}...b_{s}}^{\alpha \alpha _{2}...\alpha
_{p}a_{1}...a_{r}}+...+\vspace*{2mm}\left( \rho ,\eta \right) H_{\alpha
\gamma }^{\alpha _{p}}T_{\beta _{1}...\beta _{q}b_{1}...b_{s}}^{\alpha
_{1}...\alpha _{p-1}\alpha a_{1}...a_{r}} \\
\hspace*{8mm}-\left( \rho ,\eta \right) H_{\beta _{1}\gamma }^{\beta
}T_{\beta \beta _{2}...\beta _{q}b_{1}...b_{s}}^{\alpha _{1}...\alpha
_{p}a_{1}...a_{r}}-...-\vspace*{2mm}\left( \rho ,\eta \right) H_{\beta
_{q}\gamma }^{\beta }T_{\beta _{1}...\beta _{q-1}\beta
b_{1}...b_{s}}^{\alpha _{1}...\alpha _{p}a_{1}...a_{r}} \\
\hspace*{8mm}+\left( \rho ,\eta \right) H_{a\gamma }^{a_{1}}T_{\beta
_{1}...\beta _{q}b_{1}...b_{s}}^{\alpha _{1}...\alpha
_{p}aa_{2}...a_{r}}+...+\vspace*{2mm}\left( \rho ,\eta \right) H_{a\gamma
}^{a_{r}}T_{\beta _{1}...\beta _{q}b_{1}...b_{s}}^{\alpha _{1}...\alpha
_{p}a_{1}...a_{r-1}a} \\
\hspace*{8mm}-\left( \rho ,\eta \right) H_{b_{1}\gamma }^{b}T_{\beta
_{1}...\beta _{q}bb_{2}...b_{s}}^{\alpha _{1}...\alpha _{p}a_{1}...a_{r}}-%
\vspace*{2mm}...-\left( \rho ,\eta \right) H_{b_{s}\gamma }^{b}T_{\beta
_{1}...\beta _{q}b_{1}...b_{s-1}b}^{\alpha _{1}...\alpha _{p}a_{1}...a_{r}}%
\end{array}%
\leqno(4.6)
\end{equation*}%
\emph{and }%
\begin{equation*}
\begin{array}{l}
T_{\beta _{1}...\beta _{q}b_{1}...b_{s}}^{\alpha _{1}...\alpha
_{p}a_{1}...a_{r}}\mid _{c}=\Gamma \left( \tilde{\rho},Id_{E}\right) \left(
\overset{\cdot }{\tilde{\partial}}_{c}\right) T_{\beta _{1}...\beta
_{q}b_{1}...b_{s}}^{\alpha _{1}...\alpha _{p}a_{1}...a_{r}} \\
\hspace*{8mm}+\left( \rho ,\eta \right) V_{\alpha c}^{\alpha _{1}}T_{\beta
_{1}...\beta _{q}b_{1}...b_{s}}^{\alpha \alpha _{2}...\alpha
_{p}a_{1}...a_{r}}+...+\left( \rho ,\eta \right) V_{\alpha c}^{\alpha
_{p}}T_{\beta _{1}...\beta _{q}b_{1}...b_{s}}^{\alpha _{1}...\alpha
_{p-1}\alpha a_{1}...a_{r}}\vspace*{2mm} \\
\hspace*{8mm}-\left( \rho ,\eta \right) V_{\beta _{1}c}^{\beta }T_{\beta
\beta _{2}...\beta _{q}b_{1}...b_{s}}^{\alpha _{1}...\alpha
_{p}a_{1}...a_{r}}-...-\left( \rho ,\eta \right) V_{\beta _{q}c}^{\beta
}T_{\beta _{1}...\beta _{q-1}\beta b_{1}...b_{s}}^{\alpha _{1}...\alpha
_{p}a_{1}...a_{r}}\vspace*{2mm} \\
\hspace*{8mm}+\left( \rho ,\eta \right) V_{ac}^{a_{1}}T_{\beta _{1}...\beta
_{q}b_{1}...b_{s}}^{\alpha _{1}...\alpha _{p}aa_{2}...a_{r}}+...+\left( \rho
,\eta \right) V_{ac}^{a_{r}}T_{\beta _{1}...\beta _{q}b_{1}...b_{s}}^{\alpha
_{1}...\alpha _{p}a_{1}...a_{r-1}a} \\
\hspace*{8mm}-\left( \rho ,\eta \right) V_{b_{1}c}^{b}T_{\beta _{1}...\beta
_{q}bb_{2}...b_{s}}^{\alpha _{1}...\alpha _{p}a_{1}...a_{r}}-...-\left( \rho
,\eta \right) V_{b_{s}c}^{b}T_{\beta _{1}...\beta
_{q}b_{1}...b_{s-1}b}^{\alpha _{1}...\alpha _{p}a_{1}...a_{r}}.%
\end{array}%
\leqno(4.7)
\end{equation*}

\textbf{Corollary 4.2 }\emph{In the particular case of Lie algebroids, }$%
\left( \eta ,h\right) =\left( Id_{M},Id_{M}\right) ,$\emph{\ we obtain }%
\begin{equation*}
\begin{array}{l}
T_{\beta _{1}...\beta _{q}b_{1}...b_{s}\mid \gamma }^{\alpha _{1}...\alpha
_{p}a_{1}...a_{r}}=\vspace*{2mm}\Gamma \left( \tilde{\rho},Id_{E}\right)
\left( \tilde{\delta}_{\gamma }\right) T_{\beta _{1}...\beta
_{q}b_{1}...b_{s}}^{\alpha _{1}...\alpha _{p}a_{1}...a_{r}} \\
\hspace*{8mm}+\rho H_{\alpha \gamma }^{\alpha _{1}}T_{\beta _{1}...\beta
_{q}b_{1}...b_{s}}^{\alpha \alpha _{2}...\alpha _{p}a_{1}...a_{r}}+...+%
\vspace*{2mm}\rho H_{\alpha \gamma }^{\alpha _{p}}T_{\beta _{1}...\beta
_{q}b_{1}...b_{s}}^{\alpha _{1}...\alpha _{p-1}\alpha a_{1}...a_{r}} \\
\hspace*{8mm}-\rho H_{\beta _{1}\gamma }^{\beta }T_{\beta \beta _{2}...\beta
_{q}b_{1}...b_{s}}^{\alpha _{1}...\alpha _{p}a_{1}...a_{r}}-...-\vspace*{2mm}%
\rho H_{\beta _{q}\gamma }^{\beta }T_{\beta _{1}...\beta _{q-1}\beta
b_{1}...b_{s}}^{\alpha _{1}...\alpha _{p}a_{1}...a_{r}} \\
\hspace*{8mm}+\rho H_{a\gamma }^{a_{1}}T_{\beta _{1}...\beta
_{q}b_{1}...b_{s}}^{\alpha _{1}...\alpha _{p}aa_{2}...a_{r}}+...+\vspace*{2mm%
}\rho H_{a\gamma }^{a_{r}}T_{\beta _{1}...\beta _{q}b_{1}...b_{s}}^{\alpha
_{1}...\alpha _{p}a_{1}...a_{r-1}a} \\
\hspace*{8mm}-\rho H_{b_{1}\gamma }^{b}T_{\beta _{1}...\beta
_{q}bb_{2}...b_{s}}^{\alpha _{1}...\alpha _{p}a_{1}...a_{r}}-\vspace*{2mm}%
...-\rho H_{b_{s}\gamma }^{b}T_{\beta _{1}...\beta
_{q}b_{1}...b_{s-1}b}^{\alpha _{1}...\alpha _{p}a_{1}...a_{r}}%
\end{array}%
\leqno(4.6^{\prime })
\end{equation*}%
\emph{and }%
\begin{equation*}
\begin{array}{l}
T_{\beta _{1}...\beta _{q}b_{1}...b_{s}}^{\alpha _{1}...\alpha
_{p}a_{1}...a_{r}}\mid _{c}=\Gamma \left( \tilde{\rho},Id_{E}\right) \left(
\overset{\cdot }{\tilde{\partial}}_{c}\right) T_{\beta _{1}...\beta
_{q}b_{1}...b_{s}}^{\alpha _{1}...\alpha _{p}a_{1}...a_{r}} \\
\hspace*{8mm}+\rho V_{\alpha c}^{\alpha _{1}}T_{\beta _{1}...\beta
_{q}b_{1}...b_{s}}^{\alpha \alpha _{2}...\alpha _{p}a_{1}...a_{r}}+...+\rho
V_{\alpha c}^{\alpha _{p}}T_{\beta _{1}...\beta _{q}b_{1}...b_{s}}^{\alpha
_{1}...\alpha _{p-1}\alpha a_{1}...a_{r}} \\
\hspace*{8mm}-\rho V_{\beta _{1}c}^{\beta }T_{\beta \beta _{2}...\beta
_{q}b_{1}...b_{s}}^{\alpha _{1}...\alpha _{p}a_{1}...a_{r}}-...-\rho
V_{\beta _{q}c}^{\beta }T_{\beta _{1}...\beta _{q-1}\beta
b_{1}...b_{s}}^{\alpha _{1}...\alpha _{p}a_{1}...a_{r}} \\
\hspace*{8mm}+\rho V_{ac}^{a_{1}}T_{\beta _{1}...\beta
_{q}b_{1}...b_{s}}^{\alpha _{1}...\alpha _{p}aa_{2}...a_{r}}+...+\rho
V_{ac}^{a_{r}}T_{\beta _{1}...\beta _{q}b_{1}...b_{s}}^{\alpha _{1}...\alpha
_{p}a_{1}...a_{r-1}a}\vspace*{2mm} \\
\hspace*{8mm}-\rho V_{b_{1}c}^{b}T_{\beta _{1}...\beta
_{q}bb_{2}...b_{s}}^{\alpha _{1}...\alpha _{p}a_{1}...a_{r}}-...-\rho
V_{b_{s}c}^{b}T_{\beta _{1}...\beta _{q}b_{1}...b_{s-1}b}^{\alpha
_{1}...\alpha _{p}a_{1}...a_{r}}.%
\end{array}%
\leqno(4.7^{\prime })
\end{equation*}

\emph{In the classical case, }$\left( \rho ,\eta ,h\right) =\left(
Id_{TM},Id_{M},Id_{M}\right) ,$\emph{\ we obtain }%
\begin{equation*}
\begin{array}{l}
T_{j_{1}...j_{q}b_{1}...b_{s}\mid k}^{i_{1}...i_{p}a_{1}...a_{r}}=\vspace*{%
2mm}\delta _{k}\left(
T_{j_{1}...j_{q}b_{1}...b_{s}}^{i_{1}...i_{p}a_{1}...a_{r}}\right) \\
\hspace*{8mm}%
+H_{ik}^{i_{1}}T_{j_{1}...j_{q}b_{1}...b_{s}}^{ii_{2}...i_{p}a_{1}...a_{r}}+...+%
\vspace*{2mm}H_{ik}^{i_{p}}T_{\beta _{1}...\beta
_{q}b_{1}...b_{s}}^{i_{1}...i_{p-1}ia_{1}...a_{r}} \\
\hspace*{8mm}%
-H_{j_{1}k}^{j}T_{jj_{2}...j_{q}b_{1}...b_{s}}^{i_{1}...i_{p}a_{1}...a_{r}}-...-%
\vspace*{2mm}H_{j_{q}k}^{j}T_{j_{1}...j_{q-1}jb_{1}...b_{s}}^{\alpha
_{1}...\alpha _{p}a_{1}...a_{r}} \\
\hspace*{8mm}+H_{ak}^{a_{1}}T_{\beta _{1}...\beta _{q}b_{1}...b_{s}}^{\alpha
_{1}...\alpha _{p}aa_{2}...a_{r}}+...+\vspace*{2mm}H_{ak}^{a_{r}}T_{\beta
_{1}...\beta _{q}b_{1}...b_{s}}^{\alpha _{1}...\alpha _{p}a_{1}...a_{r-1}a}
\\
\hspace*{8mm}-H_{b_{1}k}^{b}T_{\beta _{1}...\beta
_{q}bb_{2}...b_{s}}^{\alpha _{1}...\alpha _{p}a_{1}...a_{r}}-\vspace*{2mm}%
...-H_{b_{s}k}^{b}T_{\beta _{1}...\beta _{q}b_{1}...b_{s-1}b}^{\alpha
_{1}...\alpha _{p}a_{1}...a_{r}}%
\end{array}%
\leqno(4.6^{\prime \prime })
\end{equation*}%
\emph{and }%
\begin{equation*}
\begin{array}{l}
T_{j_{1}...j_{q}b_{1}...b_{s}}^{i_{1}...i_{p}a_{1}...a_{r}}\mid _{c}=\dot{%
\partial}_{c}\left( T_{\beta _{1}...\beta _{q}b_{1}...b_{s}}^{\alpha
_{1}...\alpha _{p}a_{1}...a_{r}}\right) \vspace*{2mm} \\
\hspace*{8mm}%
+V_{ic}^{i_{1}}T_{j_{1}...j_{q}b_{1}...b_{s}}^{ii_{2}...i_{p}a_{1}...a_{r}}+...+V_{ic}^{i_{p}}T_{\beta _{1}...\beta _{q}b_{1}...b_{s}}^{i_{1}...i_{p-1}ia_{1}...a_{r}}%
\vspace*{2mm} \\
\hspace*{8mm}%
-V_{j_{1}c}^{j}T_{jj_{2}...j_{q}b_{1}...b_{s}}^{i_{1}...i_{p}a_{1}...a_{r}}-...-V_{j_{q}c}^{j}T_{j_{1}...j_{q-1}jb_{1}...b_{s}}^{i_{1}...i_{p}a_{1}...a_{r}}%
\vspace*{2mm} \\
\hspace*{8mm}%
+V_{ac}^{a_{1}}T_{j_{1}...j_{q}b_{1}...b_{s}}^{i_{1}...i_{p}aa_{2}...a_{r}}+...+V_{ac}^{a_{r}}T_{j_{1}...j_{q}b_{1}...b_{s}}^{i_{1}...i_{p}a_{1}...a_{r-1}a}%
\vspace*{2mm} \\
\hspace*{8mm}%
-V_{b_{1}c}^{b}T_{j_{1}...j_{q}bb_{2}...b_{s}}^{i_{1}...i_{p}a_{1}...a_{r}}-...-V_{b_{s}c}^{b}T_{j_{1}...j_{q}b_{1}...b_{s-1}b}^{i_{1}...i_{p}a_{1}...a_{r}}.%
\end{array}%
\leqno(4.7^{\prime \prime })
\end{equation*}

\textbf{Definition 4.2 }If $\left( E,\pi ,M\right) =\left( F,\nu ,N\right) ,$
$\left( \rho ,\eta \right) \Gamma $ is a $\left( \rho ,\eta \right) $%
-connection for the vector bundle $\left( E,\pi ,M\right) $ and
\begin{equation*}
\left( \left( \rho ,\eta \right) H_{bc}^{a},\left( \rho ,\eta \right) \tilde{%
H}_{bc}^{a},\left( \rho ,\eta \right) V_{bc}^{a},\left( \rho ,\eta \right)
\tilde{V}_{bc}^{a}\right)
\end{equation*}%
are the components of a distinguished linear $\left( \rho ,\eta \right) $%
\textit{-}connection for the generalized tangent bundle $\left( \left( \rho
,\eta \right) TE,\left( \rho ,\eta \right) \tau _{E},E\right) $ such that
\begin{equation*}
\left( \rho ,\eta \right) H_{bc}^{a}=\left( \rho ,\eta \right) \tilde{H}%
_{bc}^{a}\mbox{ and }\left( \rho ,\eta \right) V_{bc}^{a}=\left( \rho ,\eta
\right) \tilde{V}_{bc}^{a},
\end{equation*}%
then we will say that \emph{the generalized tangent bundle }$\!(\!(\rho
,\!\eta )TE,(\rho ,\!\eta )\tau _{E},\!E)$ \emph{is endowed with a normal
distinguished linear }$\left( \rho ,\eta \right) $\emph{-connection on
components }$\left( \left( \rho ,\eta \right) H_{bc}^{a},\left( \rho ,\eta
\right) V_{bc}^{a}\right) $.

In the particular case of Lie algebroids, $\left( \eta ,h\right) =\left(
Id_{M},Id_{M}\right) ,$\emph{\ }the components of a normal distinguished
linear $\left( \rho ,Id_{M}\right) $-connection $\left( \rho H,\rho V\right)
$ will be denoted $\left( \rho H_{bc}^{a},\rho V_{bc}^{a}\right) $.

In the classical case, $\left( \rho ,\eta ,h\right) =\left(
Id_{TE},Id_{M},Id_{M}\right) ,$\emph{\ }the components of a normal
distinguished linear $\left( Id_{TM},Id_{M}\right) $-connection $\left(
H,V\right) $ will be denoted $\left( H_{jk}^{i},V_{jk}^{i}\right) $.

\section{The $\left( \protect\rho ,\protect\eta \right) $%
-(pseudo)metrizability}

\ \ \

We consider the following diagram:
\begin{equation*}
\begin{array}{rcl}
E &  & \left( F,\left[ ,\right] _{F,h},\left( \rho ,\eta \right) \right) \\
\pi \downarrow &  & ~\downarrow \nu \\
M & ^{\underrightarrow{~\ \ \ \ h~\ \ \ \ }} & ~\ N%
\end{array}%
\leqno(5.1)
\end{equation*}%
where $\left( E,\pi ,M\right) \in \left\vert \mathbf{B}^{\mathbf{v}%
}\right\vert $ and $\left( \left( F,\nu ,N\right) ,\left[ ,\right]
_{F,h},\left( \rho ,\eta \right) \right) $ is a generalized Lie algebroid$.$
Let $\left( \rho ,\eta \right) \Gamma $ be a $\left( \rho ,\eta \right) $%
-connection for the vector bundle $\left( E,\pi ,M\right) $ and let $\left(
\left( \rho ,\eta \right) H,\left( \rho ,\eta \right) V\right) $ be a
distinguished linear $\left( \rho ,\eta \right) $-connection for the Lie
algebroid generalized tangent bundle%
\begin{equation*}
\begin{array}{c}
\left( \left( \left( \rho ,\eta \right) TE,\left( \rho ,\eta \right) \tau
_{E},E\right) ,\left[ ,\right] _{\left( \rho ,\eta \right) TE},\left( \tilde{%
\rho},Id_{E}\right) \right) .%
\end{array}%
\end{equation*}

\textbf{Definition 5.1} A tensor field of the type%
\begin{equation*}
G=g_{\alpha \beta }d\tilde{z}^{\alpha }\otimes d\tilde{z}^{\beta
}+g_{ab}\delta \tilde{y}^{a}\otimes \delta \tilde{y}^{b}\in \mathcal{T}%
_{22}^{00}\left( \left( \rho ,\eta \right) TE,\left( \rho ,\eta \right) \tau
_{E},E\right)
\end{equation*}%
will be called \emph{pseudometrical structure }if its components are
symmetric and the matrices $\left\Vert g_{\alpha \beta }\left( u_{x}\right)
\right\Vert $and $\left\Vert g_{ab}\left( u_{x}\right) \right\Vert $ are
nondegenerate, for any point $u_{x}\in E.$

Moreover, if the matrices $\left\Vert g_{\alpha \beta }\left( u_{x}\right)
\right\Vert $ and $\left\Vert g_{ab}\left( u_{x}\right) \right\Vert $ has
constant signature, then the tensor $d$-field $G$ will be called \emph{%
metrical structure}\textit{.}

Let
\begin{equation*}
G=g_{\alpha \beta }d\tilde{z}^{\alpha }\otimes d\tilde{z}^{\beta
}+g_{ab}\delta \tilde{y}^{a}\otimes \delta \tilde{y}^{b}
\end{equation*}%
be a (pseudo)metrical structure. If $\alpha ,\beta \in \overline{1,p}$ and $%
a,b\in \overline{1,r},$ then for any vector local $\left( m+r\right) $-chart
$\left( U,s_{U}\right) $ of $\left( E,\pi ,M\right) $, we consider the real
functions
\begin{equation*}
\begin{array}{ccc}
\pi ^{-1}\left( U\right) & ^{\underrightarrow{~\ \ \tilde{g}^{\beta \alpha
}~\ \ }} & \mathbb{R}%
\end{array}%
\end{equation*}%
and
\begin{equation*}
\begin{array}{ccc}
\pi ^{-1}\left( U\right) & ^{\underrightarrow{~\ \ \tilde{g}^{ba}~\ \ }} &
\mathbb{R}%
\end{array}%
\end{equation*}%
such that%
\begin{equation*}
\begin{array}{c}
\left\Vert \tilde{g}^{\beta \alpha }\left( u_{x}\right) \right\Vert
=\left\Vert g_{\alpha \beta }\left( u_{x}\right) \right\Vert ^{-1}%
\end{array}%
\end{equation*}%
and
\begin{equation*}
\begin{array}{c}
\left\Vert \tilde{g}^{ba}\left( u_{x}\right) \right\Vert =\left\Vert
g_{ab}\left( u_{x}\right) \right\Vert ^{-1},%
\end{array}%
\end{equation*}%
for any $u_{x}\in \pi ^{-1}\left( U\right) \backslash \left\{ 0_{x}\right\} $%
.

\textbf{Definition 5.2} If around each point $x\in M$ it exists a local
vector $m+r$-chart $\left( U,s_{U}\right) $ and a local $m$-chart $\left(
U,\xi _{U}\right) $ such that $g_{\alpha \beta }\circ s_{U}^{-1}\circ \left(
\xi _{U}^{-1}\times Id_{\mathbb{R}^{r}}\right) \left( x,y\right) $ and $%
g_{ab}\circ s_{U}^{-1}\circ \left( \xi _{U}^{-1}\times Id_{\mathbb{R}%
^{r}}\right) \left( x,y\right) $ depends only on $x$, for any $u_{x}\in \pi
^{-1}\left( U\right) ,$ then we will say that the \emph{(pseudo)metrical
structure \ }%
\begin{equation*}
G=g_{\alpha \beta }d\tilde{z}^{\alpha }\otimes d\tilde{z}^{\beta
}+g_{ab}\delta \tilde{y}^{a}\otimes \delta \tilde{y}^{b}
\end{equation*}%
\emph{is a Riemannian (pseudo)metrical structure.}\textit{\ }

If only the condition is verified:

\textquotedblright $g_{\alpha \beta }\circ s_{U}^{-1}\circ \left( \xi
_{U}^{-1}\times Id_{\mathbb{R}^{r}}\right) \left( x,y\right) $\textit{\
depends only on }$x$\textit{, for any }$u_{x}\in \pi ^{-1}\left( U\right) $%
\textit{" res\-pec\-ti\-vely \textquotedblright }$g_{ab}\circ
s_{U}^{-1}\circ \left( \xi _{U}^{-1}\times Id_{\mathbb{R}^{r}}\right) \left(
x,y\right) $\textit{\ depends only on }$x$\textit{, for any }$u_{x}\in \pi
^{-1}\left( U\right) $", then we will say that the \textit{(pseudo)metrical
structure }$G$ \textit{is a }\emph{Riemannian }$\mathcal{H}$\emph{%
-(pseudo)metrical structure}\textit{\ }respectively a \emph{Riemannian }$%
\mathcal{V}$\emph{-(pseudo)metrical structure.}

\textbf{Definition 5.3} If there exists a (pseudo)metrical structure%
\begin{equation*}
G=g_{\alpha \beta }d\tilde{z}^{\alpha }\otimes d\tilde{z}^{\beta
}+g_{ab}\delta \tilde{y}^{a}\otimes \delta \tilde{y}^{b}
\end{equation*}%
such that\emph{\ }%
\begin{equation*}
\begin{array}{c}
\left( \rho ,\eta \right) D_{X}G=0,~\forall X\in \Gamma \left( \left( \rho
,\eta \right) TE,\left( \rho ,\eta \right) \tau _{E},E\right) .%
\end{array}%
\leqno(5.2)
\end{equation*}%
then the Lie algebroid generalized tangent bundle
\begin{equation*}
\begin{array}{c}
\left( \left( \left( \rho ,\eta \right) TE,\left( \rho ,\eta \right) \tau
_{E},E\right) ,\left[ ,\right] _{\left( \rho ,\eta \right) TE},\left( \tilde{%
\rho},Id_{E}\right) \right) .%
\end{array}%
\end{equation*}%
will be called $(\rho ,\eta )$\emph{-(pseudo)metrizable}

Condition $\left( 5.2\right) $ is equivalent with the following equalities:
\begin{equation*}
\begin{array}{c}
g_{\alpha \beta \mid \gamma }=0,\,g_{ab\mid \gamma }=0,\,\,g_{\alpha \beta
}\mid _{c}=0\,,\,\,g_{ab}\mid _{c}=0.%
\end{array}%
\leqno(5.3)
\end{equation*}

If $g_{\alpha \beta \mid \gamma }{=}0$ and $\,g_{ab\mid \gamma }{=}0$, then
we will say that \emph{the Lie algebroid generalized tangent bundle }%
\begin{equation*}
\begin{array}{c}
\left( \left( \left( \rho ,\eta \right) TE,\left( \rho ,\eta \right) \tau
_{E},E\right) ,\left[ ,\right] _{\left( \rho ,\eta \right) TE},\left( \tilde{%
\rho},Id_{E}\right) \right) .%
\end{array}%
\end{equation*}
\emph{is }$\mathcal{H}$\emph{-}$(\rho ,\eta )$\emph{-(pseudo)metrizable.}

If $g_{\alpha \beta }|_{c}{=}0$ and $\,g_{ab}|_{c}{=}0$, then we will say
that \emph{the Lie algebroid generalized tangent bundle }%
\begin{equation*}
\begin{array}{c}
\left( \left( \left( \rho ,\eta \right) TE,\left( \rho ,\eta \right) \tau
_{E},E\right) ,\left[ ,\right] _{\left( \rho ,\eta \right) TE},\left( \tilde{%
\rho},Id_{E}\right) \right) .%
\end{array}%
\end{equation*}
\emph{is }$\mathcal{V}$\emph{-}$(\rho ,\eta )$\emph{-(pseudo)\-metrizable.}

\textbf{Theorem 5.1} \emph{If} $\left( \left( \rho ,\eta \right) \mathring{H}%
,\left( \rho ,\eta \right) \mathring{V}\right) $ \emph{is a distinguished
linear }$\left( \rho ,\eta \right) $\emph{-connection for the Lie algebroid
generalized tangent bundle }%
\begin{equation*}
\begin{array}{c}
\left( \left( \left( \rho ,\eta \right) TE,\left( \rho ,\eta \right) \tau
_{E},E\right) ,\left[ ,\right] _{\left( \rho ,\eta \right) TE},\left( \tilde{%
\rho},Id_{E}\right) \right) .%
\end{array}%
\end{equation*}
\emph{\ and }%
\begin{equation*}
G=g_{\alpha \beta }d\tilde{z}^{\alpha }\otimes d\tilde{z}^{\beta
}+g_{ab}\delta \tilde{y}^{a}\otimes \delta \tilde{y}^{b}
\end{equation*}
\emph{is a (pseudo)metrical structure, then the following real local
functions:}%
\begin{equation*}
(5.4)%
\begin{array}{ll}
\left( \rho ,\eta \right) H_{\beta \gamma }^{\alpha }\!\! & =\displaystyle%
\frac{1}{2}\tilde{g}^{\alpha \varepsilon }\left( \Gamma \left( \tilde{\rho}%
,Id_{E}\right) \left( \tilde{\delta}_{\gamma }\right) g_{\varepsilon \beta
}+\Gamma \left( \tilde{\rho},Id_{E}\right) \left( \tilde{\delta}_{\beta
}\right) g_{\varepsilon \gamma }-\Gamma \left( \tilde{\rho},Id_{E}\right)
\left( \tilde{\delta}_{\varepsilon }\right) g_{\beta \gamma }\right.
\vspace*{1mm} \\
& \left. +g_{\theta \varepsilon }L_{\gamma \beta }^{\theta }\circ h\circ \pi
-g_{\beta \theta }L_{\gamma \varepsilon }^{\theta }\circ h\circ \pi
-g_{\theta \gamma }L_{\beta \varepsilon }^{\theta }\circ h\circ \pi \right) ,%
\vspace*{2mm} \\
\left( \rho ,\eta \right) H_{b\gamma }^{a}\!\! & =\left( \rho ,\eta \right)
\mathring{H}_{b\gamma }^{a}+\displaystyle\frac{1}{2}\tilde{g}^{ac}g_{bc%
\overset{0}{\mid }\gamma },\vspace*{2mm} \\
\left( \rho ,\eta \right) V_{\beta c}^{\alpha }\!\! & =\left( \rho ,\eta
\right) \mathring{V}_{\beta c}^{\alpha }+\displaystyle\frac{1}{2}\tilde{g}%
^{\alpha \varepsilon }g_{\beta \varepsilon \overset{0}{\mid }c},\vspace*{2mm}
\\
\left( \rho ,\eta \right) V_{bc}^{a}\!\! & =\displaystyle\frac{1}{2}\tilde{g}%
^{ae}\left( \Gamma \left( \tilde{\rho},Id_{E}\right) \left( \overset{\cdot }{%
\tilde{\partial}}_{c}\right) g_{eb}+\Gamma \left( \tilde{\rho},Id_{E}\right)
\left( \overset{\cdot }{\tilde{\partial}}_{b}\right) g_{ec}-\Gamma \left(
\tilde{\rho},Id_{E}\right) \left( \overset{\cdot }{\tilde{\partial}}%
_{e}\right) g_{bc}\right)%
\end{array}%
\end{equation*}%
\emph{are components of a distinguished linear }$\left( \rho ,\eta \right) $%
\emph{-connection such that the Lie algebroid generalized tangent bundle}
\emph{becomes }$\left( \rho ,\eta \right) $\emph{-(pseudo)metrizable.}

\textbf{Corollary 5.1 }\emph{In the particular case of Lie algebroids, }$%
\left( \eta ,h\right) =\left( Id_{M},Id_{M}\right) ,$\emph{\ then we obtain}
\begin{equation*}
(5.4^{\prime })%
\begin{array}{ll}
\rho H_{\beta \gamma }^{\alpha }\!\! & =\displaystyle\frac{1}{2}\tilde{g}%
^{\alpha \varepsilon }\left( \Gamma \left( \tilde{\rho},Id_{E}\right) \left(
\tilde{\delta}_{\gamma }\right) g_{\varepsilon \beta }+\Gamma \left( \tilde{%
\rho},Id_{E}\right) \left( \tilde{\delta}_{\beta }\right) g_{\varepsilon
\gamma }-\Gamma \left( \tilde{\rho},Id_{E}\right) \left( \tilde{\delta}%
_{\varepsilon }\right) g_{\beta \gamma }\right. \vspace*{1mm} \\
& \left. +g_{\theta \varepsilon }L_{\gamma \beta }^{\theta }\circ \pi
-g_{\beta \theta }L_{\gamma \varepsilon }^{\theta }\circ \pi -g_{\theta
\gamma }L_{\beta \varepsilon }^{\theta }\circ \pi \right) ,\vspace*{2mm} \\
\rho H_{b\gamma }^{a}\!\! & =\rho \mathring{H}_{b\gamma }^{a}+\displaystyle%
\frac{1}{2}\tilde{g}^{ac}g_{bc\overset{0}{\mid }\gamma },\vspace*{2mm} \\
\rho V_{\beta c}^{\alpha }\!\! & =\rho \mathring{V}_{\beta c}^{\alpha }+%
\displaystyle\frac{1}{2}\tilde{g}^{\alpha \varepsilon }g_{\beta \varepsilon
\overset{0}{\mid }c},\vspace*{2mm} \\
\rho V_{bc}^{a}\!\! & =\displaystyle\frac{1}{2}\tilde{g}^{ae}\left( \Gamma
\left( \tilde{\rho},Id_{E}\right) \left( \overset{\cdot }{\tilde{\partial}}%
_{c}\right) g_{eb}+\Gamma \left( \tilde{\rho},Id_{E}\right) \left( \overset{%
\cdot }{\tilde{\partial}}_{b}\right) g_{ec}-\Gamma \left( \tilde{\rho}%
,Id_{E}\right) \left( \overset{\cdot }{\tilde{\partial}}_{e}\right)
g_{bc}\right) \vspace*{1mm}%
\end{array}%
\end{equation*}

\emph{In the classicale case, }$\left( \rho ,\eta ,h\right) =\left(
Id_{TM},Id_{M},Id_{M}\right) ,$\emph{\ then we obtain}
\begin{equation*}
\begin{array}{ll}
H_{jk}^{i}\!\! & =\displaystyle\frac{1}{2}\tilde{g}^{ih}\left( \delta
_{k}g_{hj}+\delta _{j}g_{hk}-\delta _{h}g_{jk}\vspace*{1mm}\right) \\
H_{bk}^{a}\!\! & =\mathring{H}_{bk}^{a}+\displaystyle\frac{1}{2}\tilde{g}%
^{ac}g_{bc\overset{0}{\mid }k},\vspace*{2mm} \\
V_{jc}^{i}\!\! & =\mathring{V}_{jc}^{i}+\displaystyle\frac{1}{2}\tilde{g}%
^{ih}g_{jh\overset{0}{\mid }c},\vspace*{2mm} \\
V_{bc}^{a}\!\! & =\displaystyle\frac{1}{2}\tilde{g}^{ae}\left( \dot{\partial}%
_{c}g_{eb}+\dot{\partial}_{b}g_{ec}-\dot{\partial}_{e}g_{bc}\right)%
\end{array}%
\leqno(4.4^{\prime \prime })
\end{equation*}

\bigskip \noindent \textbf{Theorem 5.2 }\emph{If the distinguished linear }$%
\left( \rho ,\eta \right) $\emph{-connection} $\left( \left( \rho ,\eta
\right) \mathring{H},\left( \rho ,\eta \right) \mathring{V}\right) $ \emph{%
coincides with the Berwald linear }$\left( \rho ,\eta \right) $\emph{%
-connection in the previous theorem, then the local real functions: }
\begin{equation*}
\begin{array}{ll}
\left( \rho ,\eta \right) \overset{c}{H}_{\beta \gamma }^{\alpha }\!\!\! & =%
\displaystyle\frac{1}{2}\tilde{g}^{\alpha \varepsilon }\left( \Gamma \left(
\tilde{\rho},Id_{E}\right) \left( \tilde{\delta}_{\gamma }\right)
g_{\varepsilon \beta }+\Gamma \left( \tilde{\rho},Id_{E}\right) \left(
\tilde{\delta}_{\beta }\right) g_{\varepsilon \gamma }\right. \vspace*{1mm}
\\
& -\Gamma \left( \tilde{\rho},Id_{E}\right) \left( \tilde{\delta}%
_{\varepsilon }\right) g_{\beta \gamma }+g_{\theta \varepsilon }L_{\gamma
\beta }^{\theta }\circ h\circ \pi \vspace*{1mm}\vspace*{2mm}\left. -g_{\beta
\theta }L_{\gamma \varepsilon }^{\theta }\circ h\circ \pi -g_{\theta \gamma
}L_{\beta \varepsilon }^{\theta }\circ h\circ \pi \right) ,\vspace*{2mm} \\
\left( \rho ,\eta \right) \overset{c}{H}_{b\gamma }^{a}\!\!\! & =%
\displaystyle\frac{\partial \left( \rho ,\eta \right) \Gamma _{\gamma }^{a}}{%
\partial y^{b}}+\frac{1}{2}\tilde{g}^{ac}g_{bc\overset{0}{\mid }\gamma },%
\vspace*{2mm} \\
\left( \rho ,\eta \right) \overset{c}{V}_{\beta c}^{\alpha }\!\! & =%
\displaystyle\frac{1}{2}\tilde{g}^{\alpha \varepsilon }\frac{\partial
g_{\beta \varepsilon }}{\partial y^{c}},\vspace*{2mm} \\
\left( \rho ,\eta \right) \overset{c}{V}_{bc}^{a}\!\!\! & =\displaystyle%
\frac{1}{2}\tilde{g}^{ae}\left( \frac{\partial g_{e\beta }}{\partial y^{c}}+%
\frac{\partial g_{ec}}{\partial y^{b}}-\frac{\partial g_{bc}}{\partial y^{e}}%
\right)%
\end{array}%
\hspace*{-6mm}\leqno(5.5)
\end{equation*}%
\emph{are the components of a distinguished linear }$\left( \rho ,\eta
\right) $\emph{-connection such that the Lie algebroid generalized tangent
bundle }%
\begin{equation*}
\begin{array}{c}
\left( \left( \left( \rho ,\eta \right) TE,\left( \rho ,\eta \right) \tau
_{E},E\right) ,\left[ ,\right] _{\left( \rho ,\eta \right) TE},\left( \tilde{%
\rho},Id_{E}\right) \right) .%
\end{array}%
\end{equation*}%
\emph{becomes }$\left( \rho ,\eta \right) $\emph{-(pseudo)metrizable.}

\emph{Moreover, if the (pseudo)metrical structure }$G$\emph{\ is }$\mathcal{H%
}$\emph{- and }$\mathcal{V}$\emph{-Rieman\-nian, then the local real
functions:}\smallskip \noindent
\begin{equation*}
\begin{array}{cl}
(\rho ,\eta )\overset{c}{H}_{\beta \gamma }^{\alpha } & {=}\frac{1}{2}\tilde{%
g}^{\alpha \varepsilon }\left( \rho _{\gamma }^{k}{\circ }h{\circ }\pi \frac{%
\partial g_{\varepsilon \beta }}{\partial x^{k}}+\rho _{\beta }^{j}{\circ }h{%
\circ }\pi \frac{\partial g_{\varepsilon \gamma }}{\partial x^{j}}-\rho
_{\varepsilon }^{e}{\circ }h{\circ }\pi \frac{\partial g_{\beta \gamma }}{%
\partial x^{e}}+\right. \\
& \left. +g_{\theta \varepsilon }L_{\gamma \beta }^{\theta }{\circ }h{\circ }%
\pi -g_{\beta \theta }L_{\gamma \varepsilon }^{\theta }{\circ }h{\circ }\pi
-g_{\theta \gamma }L_{\beta \varepsilon }^{\theta }{\circ }h{\circ }\pi
\right) ,\vspace*{1mm} \\
\left( \rho ,\eta \right) \overset{c}{H}_{b\gamma }^{a} & {=}\frac{\partial
\left( \rho ,\eta \right) \Gamma _{\gamma }^{a}}{\partial y^{b}}+\frac{1}{2}%
\tilde{g}^{ac}\left( \rho _{\gamma }^{i}{\circ }h{\circ }\pi \frac{\partial
g_{bc}}{\partial x^{i}}-\frac{\partial \left( \rho ,\eta \right) \Gamma
_{\gamma }^{e}}{\partial y^{b}}g_{ec}-\frac{\partial \left( \rho ,\eta
\right) \Gamma _{\gamma }^{e}}{\partial y^{c}}g_{eb}\right) , \\
\left( \rho ,\eta \right) \overset{c}{V}_{\beta c}^{\alpha } & =0, \\
\left( \rho ,\eta \right) \overset{c}{V}_{bc}^{a} & =0.%
\end{array}%
\leqno(5.6)
\end{equation*}%
\emph{are the components of a distinguished linear }$\left( \rho ,\eta
\right) $\emph{-connection such that the Lie algebroid generalized tangent
bundle} \emph{becomes }$\left( \rho ,\eta \right) $\emph{-(pseudo)metrizable.%
}

\textbf{Corollary 5.2 }\emph{In the particular case of Lie algebroids, }$%
\left( \eta ,h\right) =\left( Id_{M},Id_{M}\right) ,$\emph{\ then we obtain}
\emph{\ }
\begin{equation*}
\begin{array}{ll}
\rho \overset{c}{H}_{\beta \gamma }^{\alpha }\!\!\! & =\displaystyle\frac{1}{%
2}\tilde{g}^{\alpha \varepsilon }\left( \Gamma \left( \tilde{\rho}%
,Id_{E}\right) \left( \tilde{\delta}_{\gamma }\right) g_{\varepsilon \beta
}+\Gamma \left( \tilde{\rho},Id_{E}\right) \left( \tilde{\delta}_{\beta
}\right) g_{\varepsilon \gamma }\right. \vspace*{1mm} \\
& -\Gamma \left( \tilde{\rho},Id_{E}\right) \left( \tilde{\delta}%
_{\varepsilon }\right) g_{\beta \gamma }+g_{\theta \varepsilon }L_{\gamma
\beta }^{\theta }\circ \pi \vspace*{1mm}\left. -g_{\beta \theta }L_{\gamma
\varepsilon }^{\theta }\circ \pi -g_{\theta \gamma }L_{\beta \varepsilon
}^{\theta }\circ \pi \right) \\
\rho \overset{c}{H}_{b\gamma }^{a}\!\!\! & =\displaystyle\frac{\partial \rho
\Gamma _{\gamma }^{a}}{\partial y^{b}}+\frac{1}{2}\tilde{g}^{ac}g_{bc\overset%
{0}{\mid }\gamma },\vspace*{2mm} \\
\rho \overset{c}{V}_{\beta c}^{\alpha }\!\! & =\displaystyle\frac{1}{2}%
\tilde{g}^{\alpha \varepsilon }\frac{\partial g_{\beta \varepsilon }}{%
\partial y^{c}},\vspace*{2mm} \\
\rho \overset{c}{V}_{bc}^{a}\!\!\! & =\displaystyle\frac{1}{2}\tilde{g}%
^{ae}\left( \frac{\partial g_{e\beta }}{\partial y^{c}}+\frac{\partial g_{ec}%
}{\partial y^{b}}-\frac{\partial g_{bc}}{\partial y^{e}}\right)%
\end{array}%
\hspace*{-6mm}\leqno(5.5^{\prime })
\end{equation*}

\emph{If the (pseudo)metrical structure }$G$\emph{\ is }$\mathcal{H}$\emph{-
and }$\mathcal{V}$\emph{-Rieman\-nian, then }%
\begin{equation*}
\begin{array}{l}
\rho \overset{c}{H}_{\beta \gamma }^{\alpha }{=}\displaystyle\frac{1}{2}%
\tilde{g}^{\alpha \varepsilon }\left( \rho _{\gamma }^{k}{\circ }\pi \frac{%
\partial g_{\varepsilon \beta }}{\partial x^{k}}+\rho _{\beta }^{j}{\circ }%
\pi \frac{\partial g_{\varepsilon \gamma }}{\partial x^{j}}-\rho
_{\varepsilon }^{e}{\circ }\pi \frac{\partial g_{\beta \gamma }}{\partial
x^{e}}+\right. \vspace*{1mm} \\
\hfill \left. +g_{\theta \varepsilon }L_{\gamma \beta }^{\theta }{\circ }\pi
-g_{\beta \theta }L_{\gamma \varepsilon }^{\theta }{\circ }\pi -g_{\theta
\gamma }L_{\beta \varepsilon }^{\theta }{\circ }\pi \right) ,\vspace*{1mm}
\\
\rho \overset{c}{H}_{b\gamma }^{a}{=}\frac{\partial \rho \Gamma _{\gamma
}^{a}}{\partial y^{b}}+\frac{1}{2}\tilde{g}^{ac}\left( \rho _{\gamma }^{i}{%
\circ }\pi \frac{\partial g_{bc}}{\partial x^{i}}-\frac{\partial \rho \Gamma
_{\gamma }^{e}}{\partial y^{b}}g_{ec}-\frac{\partial \rho \Gamma _{\gamma
}^{e}}{\partial y^{c}}g_{eb}\right) ,\vspace*{2mm} \\
\rho \overset{c}{V}_{\beta c}^{\alpha }=0,\ \rho \overset{c}{V}_{bc}^{a}=0%
\end{array}%
\leqno(5.6^{\prime })
\end{equation*}

\emph{In the classicale case, }$\left( \rho ,\eta ,h\right) =\left(
Id_{TM},Id_{M},Id_{M}\right) ,$\emph{\ then we obtain}%
\begin{equation*}
\begin{array}{ll}
\overset{c}{H}_{jk}^{i}\!\!\! & =\displaystyle\frac{1}{2}\tilde{g}%
^{ih}\left( \delta _{k}g_{hj}+\delta _{j}g_{hk}-\delta _{h}g_{jk}\right) \\
\overset{c}{H}_{bk}^{a}\!\!\! & =\displaystyle\frac{\partial \Gamma _{k}^{a}%
}{\partial y^{b}}+\frac{1}{2}\tilde{g}^{ac}g_{bc\overset{0}{\mid }k},%
\vspace*{2mm} \\
\overset{c}{V}_{jc}^{i}\!\! & =\displaystyle\frac{1}{2}\tilde{g}^{ih}\frac{%
\partial g_{jh}}{\partial y^{c}},\vspace*{2mm} \\
\overset{c}{V}_{bc}^{a}\!\!\! & =\displaystyle\frac{1}{2}\tilde{g}%
^{ae}\left( \frac{\partial g_{e\beta }}{\partial y^{c}}+\frac{\partial g_{ec}%
}{\partial y^{b}}-\frac{\partial g_{bc}}{\partial y^{e}}\right)%
\end{array}%
\hspace*{-6mm}\leqno(5.5^{\prime \prime })
\end{equation*}

\emph{If the (pseudo)metrical structure }$G$\emph{\ is }$\mathcal{H}$\emph{-
and }$\mathcal{V}$\emph{-Rieman\-nian, then }%
\begin{equation*}
\begin{array}{l}
\overset{c}{H}_{jk}^{i}{=}\displaystyle\frac{1}{2}\tilde{g}^{ih}\left( \frac{%
\partial g_{hj}}{\partial x^{k}}+\frac{\partial g_{hk}}{\partial x^{j}}-%
\frac{\partial g_{jk}}{\partial x^{h}}\right) \\
\hfill \\
\overset{c}{H}_{bk}^{a}{=}\frac{\partial \Gamma _{k}^{a}}{\partial y^{b}}+%
\frac{1}{2}\tilde{g}^{ac}\left( \frac{\partial g_{bc}}{\partial x^{i}}-\frac{%
\partial \Gamma _{k}^{e}}{\partial y^{b}}g_{ec}-\frac{\partial \Gamma
_{k}^{e}}{\partial y^{c}}g_{eb}\right) ,\vspace*{2mm} \\
\overset{c}{V}_{jc}^{i}=0,\ \overset{c}{V}_{bc}^{a}=0%
\end{array}%
\leqno(4.6^{\prime \prime })
\end{equation*}%
\textbf{\ }

\textbf{Theorem 5.3} \emph{Let }%
\begin{equation*}
\left( \left( \rho ,\eta \right) \mathring{H},\left( \rho ,\eta \right)
\mathring{V}\right)
\end{equation*}%
\emph{be a distinguished linear }$\left( \rho ,\eta \right) $\emph{%
-connection for the Lie algebroid generalized tangent bundle }%
\begin{equation*}
\begin{array}{c}
\left( \left( \left( \rho ,\eta \right) TE,\left( \rho ,\eta \right) \tau
_{E},E\right) ,\left[ ,\right] _{\left( \rho ,\eta \right) TE},\left( \tilde{%
\rho},Id_{E}\right) \right) .%
\end{array}%
\end{equation*}%
\emph{\ and let }%
\begin{equation*}
\begin{array}{c}
G=g_{\alpha \beta }d\tilde{z}^{\alpha }\otimes d\tilde{z}^{\beta
}+g_{ab}\delta \tilde{y}^{a}\otimes \delta \tilde{y}^{b}%
\end{array}%
\end{equation*}%
\emph{be a (pseudo)metrical structure.}

\emph{Let }
\begin{equation*}
\begin{array}{ll}
O_{\beta \gamma }^{\alpha \varepsilon }=\displaystyle\frac{1}{2}\left(
\delta _{\beta }^{\alpha }\delta _{\gamma }^{\varepsilon }-g_{\beta \gamma }%
\tilde{g}^{\alpha \varepsilon }\right) , & O_{\beta \gamma }^{\ast \alpha
\varepsilon }=\displaystyle\frac{1}{2}\left( \delta _{\beta }^{\alpha
}\delta _{\gamma }^{\varepsilon }+g_{\beta \gamma }\tilde{g}^{\alpha
\varepsilon }\right) ,\vspace*{2mm} \\
O_{bc}^{ae}=\displaystyle\frac{1}{2}\left( \delta _{b}^{a}\delta
_{c}^{e}-g_{bc}\tilde{g}^{ae}\right) , & O_{bc}^{\ast ae}=\displaystyle\frac{%
1}{2}\left( \delta _{b}^{a}\delta _{c}^{e}+g_{bc}\tilde{g}^{ae}\right) ,%
\end{array}%
\leqno(5.7)
\end{equation*}%
\emph{be the Obata operators}.

\emph{If the real local functions }$X_{\beta \gamma }^{\alpha },X_{\beta
c}^{\alpha },Y_{b\gamma }^{a},Y_{bc}^{a}$ \emph{are components of tensor
fields,} \emph{then the local real functions given in the following: }%
\vspace*{-3mm}
\begin{equation*}
\begin{array}{l}
\left( \rho ,\eta \right) H_{\beta \gamma }^{\alpha }=\left( \rho ,\eta
\right) \overset{c}{H_{\beta \gamma }^{\alpha }}+O_{\gamma \eta }^{\alpha
\varepsilon }X_{\varepsilon \beta }^{\eta },\vspace*{1mm} \\
\left( \rho ,\eta \right) H_{b\gamma }^{a}=\left( \rho ,\eta \right) \overset%
{c}{H_{b\gamma }^{a}}+O_{bd}^{ae}Y_{e\gamma }^{d},\vspace*{1mm} \\
\left( \rho ,\eta \right) V_{\beta c}^{\alpha }=\left( \rho ,\eta \right)
\overset{c}{V_{\beta c}^{\alpha }}+\overset{\ast }{O}_{\beta \eta }^{\alpha
\varepsilon }X_{\varepsilon c}^{\eta },\vspace*{1mm} \\
\left( \rho ,\eta \right) V_{bc}^{a}=\left( \rho ,\eta \right) \overset{c}{%
V_{bc}^{a}}+\overset{\ast }{O}_{bd}^{ae}Y_{ec}^{d},%
\end{array}%
\leqno(5.8)
\end{equation*}%
\emph{are the components of a distinguished linear }$\left( \rho ,\eta
\right) $\emph{-connection such that the Lie algebroid generalized tangent
bundle} \emph{becomes }$\left( \rho ,\eta \right) $\emph{-(pseudo)metrizable.%
}

\textbf{Corollary 5.3 }\emph{In the particular case of Lie algebroids, }$%
\left( \eta ,h\right) =\left( Id_{M},Id_{M}\right) ,$\emph{\ then we obtain}%
\begin{equation*}
\begin{array}{l}
\rho H_{\beta \gamma }^{\alpha }=\rho \overset{c}{H_{\beta \gamma }^{\alpha }%
}+O_{\gamma \eta }^{\alpha \varepsilon }X_{\varepsilon \beta }^{\eta },%
\vspace*{1mm} \\
\rho H_{b\gamma }^{a}=\rho \overset{c}{H_{b\gamma }^{a}}+O_{bd}^{ae}Y_{e%
\gamma }^{d},\vspace*{1mm} \\
\rho V_{\beta c}^{\alpha }=\rho \overset{c}{V_{\beta c}^{\alpha }}+O_{\beta
\eta }^{\ast \alpha \varepsilon }X_{\varepsilon c}^{\eta },\vspace*{1mm} \\
\rho V_{bc}^{a}=\rho \overset{c}{V_{bc}^{a}}+O_{bd}^{\ast ae}Y_{ec}^{d},%
\end{array}%
\leqno(5.8^{\prime })
\end{equation*}

\emph{In the classicale case, }$\left( \rho ,\eta ,h\right) =\left(
Id_{TM},Id_{M},Id_{M}\right) ,$\emph{\ then we obtain }(see $\left[ 28\right]
$)%
\begin{equation*}
\begin{array}{l}
H_{jk}^{i}=\overset{c}{H_{jk}^{i}}+O_{kh}^{il}X_{lj}^{h},\vspace*{1mm} \\
H_{bk}^{a}=\overset{c}{H_{bk}^{a}}+O_{bd}^{ae}Y_{ek}^{d},\vspace*{1mm} \\
V_{jc}^{i}=\overset{c}{V_{jc}^{i}}+\overset{\ast }{O}_{jh}^{il}X_{lc}^{h},%
\vspace*{1mm} \\
\rho V_{bc}^{a}=\overset{c}{V_{bc}^{a}}+\overset{\ast }{O}%
_{bd}^{ae}Y_{ec}^{d},%
\end{array}%
\leqno(5.8^{\prime \prime })
\end{equation*}

\textbf{Theorem 5.4 }\emph{If }%
\begin{equation*}
\left( \left( \rho ,\eta \right) \mathring{H},\left( \rho ,\eta \right)
\mathring{V}\right)
\end{equation*}%
\emph{is a distinguished linear }$\left( \rho ,\eta \right) $\emph{%
-connection for the generalized tangent bundle }%
\begin{equation*}
\left( \left( \rho ,\eta \right) TE,\left( \rho ,\eta \right) \tau
_{E},E\right)
\end{equation*}%
\emph{\ and }%
\begin{equation*}
G=g_{\alpha \beta }d\tilde{z}^{\alpha }\otimes d\tilde{z}^{\beta
}+g_{ab}\delta \tilde{y}^{a}\otimes \delta \tilde{y}^{b}
\end{equation*}%
\emph{is a (pseudo)metrical structure, then the real local functions: }%
\begin{equation*}
\begin{array}{l}
\left( \rho ,\eta \right) H_{\beta \gamma }^{\alpha }=\left( \rho ,\eta
\right) \mathring{H}_{\beta \gamma }^{\alpha }+\displaystyle\frac{1}{2}%
\tilde{g}^{\alpha \varepsilon }g_{\varepsilon \beta \overset{0}{\mid }\gamma
},\vspace*{2mm} \\
\left( \rho ,\eta \right) H_{b\gamma }^{a}=\left( \rho ,\eta \right)
\mathring{H}_{b\gamma }^{a}+\displaystyle\frac{1}{2}\tilde{g}^{ae}g_{eb%
\overset{0}{\mid }\gamma },\vspace*{2mm} \\
\left( \rho ,\eta \right) V_{\beta c}^{\alpha }=\left( \rho ,\eta \right)
\mathring{V}_{\beta c}^{\alpha }+\displaystyle\frac{1}{2}\tilde{g}^{\alpha
\varepsilon }g_{\varepsilon \beta }\overset{0}{\mid }_{c},\vspace*{2mm} \\
\left( \rho ,\eta \right) V_{bc}^{a}=\left( \rho ,\eta \right) \mathring{V}%
_{bc}^{a}+\displaystyle\frac{1}{2}\tilde{g}^{ae}g_{eb}\overset{0}{\mid }_{c}%
\end{array}%
\leqno(5.9)
\end{equation*}%
\emph{are the components of a distinguished linear }$\left( \rho ,\eta
\right) $\emph{-connection such that the generalized tangent bundle }$\left(
\left( \rho ,\eta \right) TE,\left( \rho ,\eta \right) \tau _{E},E\right) $
\emph{becomes }$\left( \rho ,\eta \right) $\emph{-(pseudo)metrizable.}

\textbf{Corollary 5.4 }\emph{In the particular case of Lie algebroids, }$%
\left( \eta ,h\right) =\left( Id_{M},Id_{M}\right) ,$\emph{\ then we obtain}%
\begin{equation*}
\begin{array}{l}
\rho H_{\beta \gamma }^{\alpha }=\rho \mathring{H}_{\beta \gamma }^{\alpha }+%
\displaystyle\frac{1}{2}\tilde{g}^{\alpha \varepsilon }g_{\varepsilon \beta
\overset{0}{\mid }\gamma },\vspace*{2mm} \\
\rho H_{b\gamma }^{a}=\rho \mathring{H}_{b\gamma }^{a}+\displaystyle\frac{1}{%
2}\tilde{g}^{ae}g_{eb\overset{0}{\mid }\gamma },\vspace*{2mm} \\
\rho V_{\beta c}^{\alpha }=\rho \mathring{V}_{\beta c}^{\alpha }+%
\displaystyle\frac{1}{2}\tilde{g}^{\alpha \varepsilon }g_{\varepsilon \beta }%
\overset{0}{\mid }_{c},\vspace*{2mm} \\
\rho V_{bc}^{a}=\rho \mathring{V}_{bc}^{a}+\displaystyle\frac{1}{2}\tilde{g}%
^{ae}g_{eb}\overset{0}{\mid }_{c}%
\end{array}%
\leqno(5.9^{\prime })
\end{equation*}

\emph{In the classicale case, }$\left( \rho ,\eta ,h\right) =\left(
Id_{TM},Id_{M},Id_{M}\right) ,$\emph{\ then we obtain }(see $\left[ 28\right]
$)%
\begin{equation*}
\begin{array}{l}
H_{jk}^{i}=\mathring{H}_{jk}^{i}+\displaystyle\frac{1}{2}\tilde{g}^{ih}g_{hj%
\overset{0}{\mid }k},\vspace*{2mm} \\
H_{bk}^{a}=\mathring{H}_{bk}^{a}+\displaystyle\frac{1}{2}\tilde{g}^{ae}g_{eb%
\overset{0}{\mid }k},\vspace*{2mm} \\
V_{jc}^{i}=\mathring{V}_{jc}^{i}+\displaystyle\frac{1}{2}\tilde{g}^{ih}g_{hj}%
\overset{0}{\mid }_{c},\vspace*{2mm} \\
V_{bc}^{a}=\mathring{V}_{bc}^{a}+\displaystyle\frac{1}{2}\tilde{g}^{ae}g_{eb}%
\overset{0}{\mid }_{c}%
\end{array}%
\leqno(5.9^{\prime \prime })
\end{equation*}

\section{(Generalized) Lagrange $\left( \protect\rho ,\protect\eta \right) $%
-spaces}

\ \ \

We consider the following diagram:
\begin{equation*}
\begin{array}{rcl}
E &  & \left( E,\left[ ,\right] _{E,h},\left( \rho ,\eta \right) \right) \\
\pi \downarrow &  & ~\downarrow \pi \\
M & ^{\underrightarrow{~\ \ \ \ h~\ \ \ \ }} & ~\ M%
\end{array}%
\leqno(6.1)
\end{equation*}%
where $\left( \left( E,\pi ,M\right) ,\left[ ,\right] _{E,h},\left( \rho
,\eta \right) \right) $ is a generalized Lie algebroid. Let $\left( \rho
,\eta \right) \Gamma $ be a $\left( \rho ,\eta \right) $-connection for the
vector bundle $\left( E,\pi ,M\right) .$

We admit that the Lie algebroid generalized tangent bundle%
\begin{equation*}
\begin{array}{c}
\left( \left( \left( \rho ,\eta \right) TE,\left( \rho ,\eta \right) \tau
_{E},E\right) ,\left[ ,\right] _{\left( \rho ,\eta \right) TE},\left( \tilde{%
\rho},Id_{E}\right) \right)%
\end{array}%
,
\end{equation*}%
is $\left( \rho ,\eta \right) $-(pseudo)metrizable.

Let%
\begin{equation*}
G=g_{ab}d\tilde{z}^{a}\otimes d\tilde{z}^{a}+g_{ab}\delta \tilde{y}%
^{a}\otimes \delta \tilde{y}^{b}
\end{equation*}%
be a (pseudo)metrical structure and let
\begin{equation*}
\left( \left( \rho ,\eta \right) H,\left( \rho ,\eta \right) V\right)
\end{equation*}%
be a distinguished linear $\left( \rho ,\eta \right) $-connection such that%
\emph{\ }%
\begin{equation*}
\begin{array}{c}
\left( \rho ,\eta \right) D_{X}G=0,~\forall X\in \Gamma \left( \left( \rho
,\eta \right) TE,\left( \rho ,\eta \right) \tau _{E},E\right) .%
\end{array}%
\end{equation*}

\textbf{Definition 6.1 }If the (pseudo)metrical structure $G$\ is determined
by the help of a (pseudo)metrical structure
\begin{equation*}
\begin{array}{c}
g\in \mathcal{T}~_{2}^{0}\left( V\left( \rho ,\eta \right) TE,\left( \rho
,\eta \right) ,{\tau _{E}},E\right) ,%
\end{array}%
\end{equation*}%
then the $\left( \rho ,\eta \right) $-(pseudo)metrizable Lie algebroid
generalized tangent bundle%
\begin{equation*}
\begin{array}{c}
\left( \left( \left( \rho ,\eta \right) TE,\left( \rho ,\eta \right) \tau
_{E},E\right) ,\left[ ,\right] _{\left( \rho ,\eta \right) TE},\left( \tilde{%
\rho},Id_{E}\right) \right)%
\end{array}%
\end{equation*}%
will be called the \emph{generalized Lagrange }$\left( \rho ,\eta \right) $%
\emph{-space.}

\emph{Remark 6.1 The generalized Lagrange }$\left( Id_{TM},Id_{M}\right) $%
\emph{-spaces are the usual generalized Lagrange spaces.}

\textbf{Theorem 6.1 }\emph{If the (pseudo)metrical structure }$G$\emph{\ is
determined by a (pseudo)metrical structure }%
\begin{equation*}
\begin{array}{c}
g\in \mathcal{T}~_{2}^{0}\left( V\left( \rho ,\eta \right) TE,\left( \rho
,\eta \right) ,{\tau _{E}},E\right) ,%
\end{array}%
\end{equation*}%
\emph{then, the real local functions:}
\begin{equation*}
(6.1)%
\begin{array}{ll}
\left( \rho ,\eta \right) H_{bc}^{a}\!\!\! & =\displaystyle\frac{1}{2}\tilde{%
g}^{ae}\left( \Gamma \left( \tilde{\rho},Id_{E}\right) \left( \tilde{\delta}%
_{b}\right) g_{ec}+\Gamma \left( \tilde{\rho},Id_{E}\right) \left( \tilde{%
\delta}_{c}\right) g_{be}\right. -\Gamma \left( \tilde{\rho},Id_{E}\right)
\left( \tilde{\delta}_{e}\right) g_{bc}\vspace*{1mm} \\
& -\,g_{cd}L_{be}^{d}{\circ }h{\circ }\pi \left. +g_{bd}L_{ec}^{d}{\circ }h{%
\circ }\pi -g_{ed}L_{bc}^{d}{\circ }h{\circ }\pi \right) ,\vspace*{2mm} \\
\left( \rho ,\eta \right) V_{bc}^{a}\!\!\! & =\displaystyle\frac{1}{2}\tilde{%
g}^{ae}\left( \Gamma \left( \tilde{\rho},Id_{E}\right) \left( \overset{\cdot
}{\tilde{\partial}}_{c}\right) g_{eb}+\Gamma \left( \tilde{\rho}%
,Id_{E}\right) \left( \overset{\cdot }{\tilde{\partial}}_{b}\right)
g_{ec}-\Gamma \left( \tilde{\rho},Id_{E}\right) \left( \overset{\cdot }{%
\tilde{\partial}}_{e}\right) g_{bc}\right)%
\end{array}%
\end{equation*}%
\medskip \emph{are the components of a normal distinguished linear }$\left(
\rho ,\eta \right) $\emph{-connection with }$\left( \rho ,\eta \right) $%
\emph{-}$\mathcal{H}\left( \mathcal{HH}\right) $\emph{\ and }$\left( \rho
,\eta \right) $\emph{-}$\mathcal{V}\left( \mathcal{VV}\right) $\emph{\
torsions free such that the Lie algebroid generalized tangent bundle }%
\begin{equation*}
\begin{array}{c}
\left( \left( \left( \rho ,\eta \right) TE,\left( \rho ,\eta \right) \tau
_{E},E\right) ,\left[ ,\right] _{\left( \rho ,\eta \right) TE},\left( \tilde{%
\rho},Id_{E}\right) \right)%
\end{array}%
\end{equation*}%
\emph{\ becomes generalized Lagrange }$\left( \rho ,\eta \right) $\emph{%
-space.}\medskip

This normal distinguished linear $(\rho ,\eta )$-connection will be called
\emph{generalized linear} $(\rho ,\eta )$\emph{-connection of Levi-Civita
type. }

\textbf{Corolary 6.1 }\emph{In the particular case of Lie algebroids, }$%
\left( \eta ,h\right) =\left( Id_{M},Id_{M}\right) ,$\emph{\ then we obtain}
\begin{equation*}
(6.1^{\prime })%
\begin{array}{ll}
\rho H_{bc}^{a}\!\!\! & =\displaystyle\frac{1}{2}\tilde{g}^{ae}\left( \Gamma
\left( \tilde{\rho},Id_{E}\right) \left( \tilde{\delta}_{b}\right)
g_{ec}+\Gamma \left( \tilde{\rho},Id_{E}\right) \left( \tilde{\delta}%
_{c}\right) g_{be}\right. -\Gamma \left( \tilde{\rho},Id_{E}\right) \left(
\tilde{\delta}_{e}\right) g_{bc}\vspace*{1mm} \\
& -\,g_{cd}L_{be}^{d}{\circ }\pi \left. +g_{bd}L_{ec}^{d}{\circ }\pi
-g_{ed}L_{bc}^{d}{\circ }\pi \right) ,\vspace*{2mm} \\
\rho V_{bc}^{a}\!\!\! & =\displaystyle\frac{1}{2}\tilde{g}^{ae}\left( \Gamma
\left( \tilde{\rho},Id_{E}\right) \left( \overset{\cdot }{\tilde{\partial}}%
_{c}\right) g_{eb}+\Gamma \left( \tilde{\rho},Id_{E}\right) \left( \overset{%
\cdot }{\tilde{\partial}}_{b}\right) g_{ec}-\Gamma \left( \tilde{\rho}%
,Id_{E}\right) \left( \overset{\cdot }{\tilde{\partial}}_{e}\right)
g_{bc}\right)%
\end{array}%
\end{equation*}

\emph{In the classicale case, }$\left( \rho ,\eta ,h\right) =\left(
Id_{TM},Id_{M},Id_{M}\right) ,$\emph{\ then we obtain}
\begin{equation*}
\begin{array}{ll}
H_{bc}^{a}\!\!\! & =\frac{1}{2}\tilde{g}^{ae}\left( \delta _{b}g_{ec}+\delta
_{c}g_{be}-\delta _{e}g_{bc}\vspace*{1mm}\right) \\
V_{bc}^{a}\!\!\! & =\displaystyle\frac{1}{2}\tilde{g}^{ae}\left( \dot{%
\partial}_{c}g_{eb}+\dot{\partial}_{b}g_{ec}-\dot{\partial}_{e}g_{bc}\right)%
\end{array}%
\leqno(6.1^{\prime \prime })
\end{equation*}

\emph{Moreover, if }$\left( E,\pi ,M\right) =\left( TM,\tau _{M},M\right) ,$%
\emph{\ then we obtain}%
\begin{equation*}
\begin{array}{ll}
H_{jk}^{i}\!\!\! & =\frac{1}{2}\tilde{g}^{ih}\left( \delta _{j}g_{hk}+\delta
_{k}g_{jh}-\delta _{h}g_{jk}\vspace*{1mm}\right) \\
V_{jk}^{i}\!\!\! & =\displaystyle\frac{1}{2}\tilde{g}^{ih}\left( \dot{%
\partial}_{k}g_{hj}+\dot{\partial}_{j}g_{hk}-\dot{\partial}_{h}g_{jk}\right)%
\end{array}%
\leqno(6.1^{\prime \prime \prime })
\end{equation*}

\textbf{Theorem 6.2 }\emph{Let\ }$\left( \left( \rho ,\eta \right) H,\left(
\rho ,\eta \right) V\right) $\emph{\ be the normal dis\-tin\-guished linear }%
$\left( \rho ,\eta \right) $\emph{-connec\-tion presented in the previous
theorem. }

\emph{If }%
\begin{equation*}
\mathbb{T}_{bc}^{a}\tilde{\delta}_{a}\otimes d\tilde{z}^{b}\otimes d\tilde{z}%
^{c}\in \mathcal{T}_{20}^{10}\left( \left( \rho ,\eta \right) TE,\left( \rho
,\eta \right) \tau _{E},E\right)
\end{equation*}%
\emph{and }%
\begin{equation*}
\mathbb{S}_{bc}^{a}\overset{\cdot }{\tilde{\partial}}_{a}\otimes \delta
\tilde{y}^{b}\otimes \delta \tilde{y}^{c}\in \mathcal{T}_{02}^{01}\left(
\left( \rho ,\eta \right) TE,\left( \rho ,\eta \right) \tau _{E},E\right)
\end{equation*}%
\emph{such that they satisfy the conditions:}%
\begin{equation*}
\mathbb{T}_{bc}^{a}=-\mathbb{T}_{cb}^{a}\wedge ~\mathbb{S}_{bc}^{a}=-\mathbb{%
S}_{cb}^{a},~\forall b,c\in \overline{1,n},
\end{equation*}%
\emph{then the following real local functions:\ }%
\begin{equation*}
\begin{array}{l}
\left( \rho ,\eta \right) \tilde{H}_{bc}^{a}=\left( \rho ,\eta \right)
H_{bc}^{a}+\displaystyle\frac{1}{2}\tilde{g}^{ae}\left( g_{ed}\mathbb{T}%
_{bc}^{d}-g_{bd}\mathbb{T}_{ec}^{d}+g_{cd}\mathbb{T}_{be}^{d}\right) ,%
\vspace*{2mm} \\
\left( \rho ,\eta \right) \tilde{V}_{bc}^{a}=\left( \rho ,\eta \right)
V_{bc}^{a}+\displaystyle\frac{1}{2}\tilde{g}^{ae}\left( g_{ed}\mathbb{S}%
_{bc}^{d}-g_{bd}\mathbb{S}_{ec}^{d}+g_{cd}\mathbb{S}_{be}^{d}\right)%
\end{array}%
\leqno(6.2)
\end{equation*}%
\emph{are the components of a normal distinguished linear }$\left( \rho
,\eta \right) $\emph{-connection with }$\left( \rho ,\eta \right) $\emph{-}$%
\mathcal{H}\left( \mathcal{HH}\right) $\emph{\ and }$\left( \rho ,\eta
\right) $\emph{-}$\mathcal{V}\left( \mathcal{VV}\right) $\emph{\ torsions a
priori given such that the Lie algebroid generalized tangent bundle }%
\begin{equation*}
\begin{array}{c}
\left( \left( \left( \rho ,\eta \right) TE,\left( \rho ,\eta \right) \tau
_{E},E\right) ,\left[ ,\right] _{\left( \rho ,\eta \right) TE},\left( \tilde{%
\rho},Id_{E}\right) \right)%
\end{array}%
\end{equation*}%
\emph{\ becomes generalized Lagrange }$\left( \rho ,\eta \right) $\emph{%
-space.}

\emph{Moreover, we obtain: }%
\begin{equation*}
\begin{array}{l}
\mathbb{T}_{bc}^{a}=\left( \rho ,\eta \right) \tilde{H}_{bc}^{a}-\left( \rho
,\eta \right) \tilde{H}_{cb}^{a}-L_{bc}^{a}\circ h\circ \pi ,\vspace*{2mm}
\\
\mathbb{S}_{bc}^{a}=\left( \rho ,\eta \right) \tilde{V}_{bc}^{a}-\left( \rho
,\eta \right) \tilde{V}_{cb}^{a}.%
\end{array}%
\leqno\left( 6.3\right)
\end{equation*}

\textbf{Corollary 6.2} \emph{In the particular case of Lie algebroids, }$%
\left( \eta ,h\right) =\left( Id_{M},Id_{M}\right) ,$\emph{\ then we obtain}%
\begin{equation*}
\begin{array}{l}
\rho \tilde{H}_{bc}^{a}=\rho H_{bc}^{a}+\displaystyle\frac{1}{2}\tilde{g}%
^{ae}\left( g_{ed}\mathbb{T}_{bc}^{d}-g_{bd}\mathbb{T}_{ec}^{d}+g_{cd}%
\mathbb{T}_{be}^{d}\right) ,\vspace*{2mm} \\
\rho \tilde{V}_{bc}^{a}=\rho V_{bc}^{a}+\displaystyle\frac{1}{2}\tilde{g}%
^{ae}\left( g_{ed}\mathbb{S}_{bc}^{d}-g_{bd}\mathbb{S}_{ec}^{d}+g_{cd}%
\mathbb{S}_{be}^{d}\right)%
\end{array}%
\leqno(6.2^{\prime })
\end{equation*}%
\emph{and }%
\begin{equation*}
\begin{array}{l}
\mathbb{T}_{bc}^{a}=\rho \tilde{H}_{bc}^{a}-\rho \tilde{H}%
_{cb}^{a}-L_{bc}^{a}\circ \pi ,\vspace*{2mm} \\
\mathbb{S}_{bc}^{a}=\rho \tilde{V}_{bc}^{a}-\rho \tilde{V}_{cb}^{a}.%
\end{array}%
\leqno\left( 6.3^{\prime }\right)
\end{equation*}

\emph{In the classicale case, }$\left( \rho ,\eta ,h\right) =\left(
Id_{TM},Id_{M},Id_{M}\right) ,$\emph{\ then we obtain}
\begin{equation*}
\begin{array}{l}
\tilde{H}_{bc}^{a}=H_{bc}^{a}+\displaystyle\frac{1}{2}\tilde{g}^{ae}\left(
g_{ed}\mathbb{T}_{bc}^{d}-g_{bd}\mathbb{T}_{ec}^{d}+g_{cd}\mathbb{T}%
_{be}^{d}\right) ,\vspace*{2mm} \\
\tilde{V}_{bc}^{a}=V_{bc}^{a}+\displaystyle\frac{1}{2}\tilde{g}^{ae}\left(
g_{ed}\mathbb{S}_{bc}^{d}-g_{bd}\mathbb{S}_{ec}^{d}+g_{cd}\mathbb{S}%
_{be}^{d}\right)%
\end{array}%
\leqno(6.2^{\prime \prime })
\end{equation*}%
\emph{and }%
\begin{equation*}
\begin{array}{l}
\mathbb{T}_{bc}^{a}=\tilde{H}_{bc}^{a}-\tilde{H}_{cb}^{a},\vspace*{2mm} \\
\mathbb{S}_{bc}^{a}=\tilde{V}_{bc}^{a}-\tilde{V}_{cb}^{a}.%
\end{array}%
\leqno\left( 6.3^{\prime \prime }\right)
\end{equation*}

\emph{In particular, if }$\left( E,\pi ,M\right) =\left( TM,\tau
_{M},M\right) ,$\emph{\ then we obtain}
\begin{equation*}
\begin{array}{l}
\tilde{H}_{jk}^{i}=H_{jk}^{i}+\displaystyle\frac{1}{2}\tilde{g}^{ie}\left(
g_{eh}\mathbb{T}_{jk}^{h}-g_{jh}\mathbb{T}_{ek}^{h}+g_{kh}\mathbb{T}%
_{je}^{h}\right) ,\vspace*{2mm} \\
\tilde{V}_{jk}^{i}=V_{jk}^{i}+\displaystyle\frac{1}{2}\tilde{g}^{ie}\left(
g_{eh}\mathbb{S}_{jk}^{h}-g_{jh}\mathbb{S}_{ek}^{h}+g_{kh}\mathbb{S}%
_{je}^{h}\right)%
\end{array}%
\leqno(6.2^{\prime \prime \prime })
\end{equation*}%
\emph{and }%
\begin{equation*}
\begin{array}{l}
\mathbb{T}_{jk}^{i}=\tilde{H}_{jk}^{i}-\tilde{H}_{kj}^{i},\vspace*{2mm} \\
\mathbb{S}_{jk}^{i}=\tilde{V}_{jk}^{i}-\tilde{V}_{kj}^{i}.%
\end{array}%
\leqno\left( 6.3^{\prime \prime \prime }\right)
\end{equation*}

We know that in the classical theory, the morphisms used for the standard
Lie algebroid
\begin{equation*}
\begin{array}{c}
\left( \left( TM,\tau _{M},M\right) ,\left[ ,\right] _{TM},\left(
Id_{TM},Id_{M}\right) \right)%
\end{array}%
\end{equation*}%
are identities. So, the Lie algebroid structure was not explicitely used.

As, for the generalized Lie algebroid%
\begin{equation*}
\begin{array}{c}
\left( \left( E,\pi ,M\right) ,\left[ ,\right] _{E,h},\left( \rho ,\eta
\right) \right)%
\end{array}%
\end{equation*}%
the morphisms used are different by the identities morphisms, it is natural
to extend the notion of Lagrange (Finsler) fundamental function from the
standard Lie algebroid
\begin{equation*}
\begin{array}{c}
\left( \left( TM,\tau _{M},M\right) ,\left[ ,\right] _{TM},\left(
Id_{TM},Id_{M}\right) \right)%
\end{array}%
\end{equation*}%
to the generalized Lie algebroid%
\begin{equation*}
\begin{array}{c}
\left( \left( E,\pi ,M\right) ,\left[ ,\right] _{E,h},\left( \rho ,\eta
\right) \right) .%
\end{array}%
\end{equation*}

\textbf{Definition 6.2 }A \emph{smooth Lagrange fundamental function} on the
generalized Lie algebroid
\begin{equation*}
\begin{array}{c}
\left( \left( E,\pi ,M\right) ,\left[ ,\right] _{E,h},\left( \rho ,\eta
\right) \right)
\end{array}%
\end{equation*}%
is a mapping $E~\ ^{\underrightarrow{\ \ L\ \ }}~\ \mathbb{R}$ which
satisfies the following conditions:\medskip

1. $L\circ u\in C^{\infty }\left( M\right) $, for any $u\in \Gamma \left(
E,\pi ,M\right) \setminus \left\{ 0\right\} $;\smallskip

2. $L\circ 0\in C^{0}\left( M\right) $, where $0$ means the null section of $%
\left( E,\pi ,M\right) .$

If, for any local vector $m+r$-chart $\left( U,s_{U}\right) $ of $\left(
E,\pi ,M\right) ,$ we have:
\begin{equation*}
\begin{array}{c}
rank\left\Vert L_{ab}\left( u_{x}\right) \right\Vert =r,%
\end{array}%
\leqno(6.4)
\end{equation*}%
for any $u_{x}\in \pi ^{-1}\left( U\right) \backslash \left\{ 0_{x}\right\} $%
, then we say that \emph{the Lagrangian }$L$\emph{\ is regular.}

\textbf{Definition 6.3 }A \emph{smooth Finsler fundamental function} on the
generalized Lie algebroid
\begin{equation*}
\begin{array}{c}
\left( \left( E,\pi ,M\right) ,\left[ ,\right] _{E,h},\left( \rho ,\eta
\right) \right)
\end{array}%
\end{equation*}%
is a mapping $E~\ ^{\underrightarrow{\ \ F\ \ }}~\ \mathbb{R}_{+}$ which
satisfies the following conditions:\medskip

1. $F\circ u\in C^{\infty }\left( M\right) $, for any $u\in \Gamma \left(
E,\pi ,M\right) \setminus \left\{ 0\right\} $;\smallskip

2. $F\circ 0\in C^{0}\left( M\right) $, where $0$ means the null section of $%
\left( E,\pi ,M\right) $;\smallskip

3. $F$ is positively $1$-homogenous on the fibres of the vector bundle $%
\left( E,\pi ,M\right) ;$\smallskip

4. For any vector local $m+r$-chart $\left( U,s_{U}\right) $ of $\left(
E,\pi ,M\right) ,$ the hessian:%
\begin{equation*}
\begin{array}{c}
\left\Vert \frac{\partial ^{2}F^{2}\left( u_{x}\right) }{\partial
y^{a}\partial y^{b}}\right\Vert%
\end{array}%
\leqno(6.5)
\end{equation*}%
is positively define for any $u_{x}\in \pi ^{-1}\left( U\right) \backslash
\left\{ 0_{x}\right\} $.

\textbf{Definition 6.4 }If the (pseudo)metrical structure $G$\ of the \emph{%
generalized Lagrange }$\left( \rho ,\eta \right) $\emph{-space }%
\begin{equation*}
\begin{array}{c}
\left( \left( \left( \rho ,\eta \right) TE,\left( \rho ,\eta \right) \tau
_{E},E\right) ,\left[ ,\right] _{\left( \rho ,\eta \right) TE},\left( \tilde{%
\rho},Id_{E}\right) \right)%
\end{array}%
\end{equation*}%
is determined by the help of the metrical structure
\begin{equation*}
\begin{array}{c}
g=\frac{1}{2}\cdot \frac{\partial ^{2}L}{\partial y^{a}\partial y^{b}}%
dy^{a}\otimes dy^{b}\in \mathcal{T}~_{2}^{0}\left( V\left( \rho ,\eta
\right) TE,\left( \rho ,\eta \right) ,{\tau _{E}},E\right) ,%
\end{array}%
\end{equation*}%
respectively
\begin{equation*}
\begin{array}{c}
g=\frac{1}{2}\cdot \frac{\partial ^{2}F^{2}}{\partial y^{a}\partial y^{b}}%
dy^{a}\otimes dy^{b}\in \mathcal{T}~_{2}^{0}\left( V\left( \rho ,\eta
\right) TE,\left( \rho ,\eta \right) ,{\tau _{E}},E\right) ,%
\end{array}%
\end{equation*}%
then the Lie algebroid generalized tangent bundle
\begin{equation*}
\begin{array}{c}
\left( \left( \left( \rho ,\eta \right) TE,\left( \rho ,\eta \right) \tau
_{E},E\right) ,\left[ ,\right] _{\left( \rho ,\eta \right) TE},\left( \tilde{%
\rho},Id_{E}\right) \right)%
\end{array}%
\end{equation*}%
will be called the \emph{Lagrange }respectively\emph{\ Finsler }$\left( \rho
,\eta \right) $\emph{-space.}

\emph{Remark 6.2} \emph{The Lagrange }$\left( Id_{TM},Id_{M}\right) $\emph{%
-spaces and the Finsler }$\left( Id_{TM},Id_{M}\right) $\emph{-spaces are
the usual Lagrange spaces and Finsler spaces.}

\ \ \emph{\ }

\section*{Acknowledgment}

\addcontentsline{toc}{section}{Acknowledgment}

I would like to thank R\u{a}dine\c{s}ti-Gorj Cultural Scientifique Society
for financial support. In memory of Prof. Dr. Gheorghe RADU and Acad. Dr.
Doc. Cornelius RADU. Dedicated to Acad. Dr. Doc. Radu MIRON from Iasy
University, Romania at his $86$-th anniversary.

\hfill

\end{document}